%
%
\input amstex.tex
\documentstyle{amsppt}
\magnification=1200
\baselineskip=13pt
\hsize=6.5truein
\vsize=8.9truein
\def\R{{\Bbb R}}
\def\N{{\Bbb N}}
\def\C{{\Bbb C}}
\def\Z{{\Bbb Z}}
\def\Zp{{\Bbb Z}_+}
\def\L{{\Cal L}}
\def\vp{\varphi}
\def\x{\left( \frac{x+x^{-1}}2 \right)}
\def\xs{\bigl( (x+x^{-1})/2 \bigr)}
\topmatter
\title Orthogonal polynomials and Laurent polynomials related
to the Hahn-Exton $q$-Bessel function
\endtitle
\rightheadtext{Orthogonal polynomials and Laurent polynomials}
\author H.T. Koelink and W. Van Assche\endauthor
\affil  Katholieke Universiteit Leuven\endaffil
\address Departement Wiskunde, Katholieke Universiteit Leuven,
Celestijnenlaan 200 B, B-3001 Leuven (Heverlee), Belgium\endaddress
\email erik\%twi\%wis\@cc3.KULeuven.ac.be,\
walter\%twi\%wis\@cc3.KULeuven.ac.be\endemail
\date May 4, 1994\enddate
\thanks The first author is supported by a NATO-Science
Fellowship of the Netherlands
Organization for Scientific Research (NWO).
The second author is a Senior Research Associate of the Belgian
National Fund for Scientific Research (NFWO).
\endthanks
\keywords orthogonal polynomials, orthogonal Laurent polynomials,
recurrence relation, Hahn-Exton $q$-Bessel function, $q$-Lommel
polynomials, Al-Salam--Chihara polynomials, continuous $q$-Hermite
polynomials, Chebyshev polynomials, perturbation
\endkeywords
\subjclass  42C05, 33D15
\endsubjclass
\abstract Laurent polynomials related to the Hahn-Exton $q$-Bessel
function, which are $q$-analogues of the Lommel polynomials,
have been introduced by Koelink and Swarttouw. The explicit strong moment
functional with respect to which the Laurent $q$-Lommel polynomials are
orthogonal is given. The strong moment functional gives rise to two positive
definite moment functionals. For the corresponding sets of orthogonal
polynomials the orthogonality measure is determined using the three-term
recurrence relation as a starting point. The relation between Chebyshev
polynomials of the second kind and the Laurent $q$-Lommel polynomials and
related functions is used to obtain estimates for the latter.
\endabstract
\endtopmatter
\document


\head 1. Introduction and motivation\endhead

The Lommel polynomials are orthogonal polynomials closely related to the
Bessel function. Although the Lommel polynomials have a representation
involving a hypergeometric ${}_2F_3$-series, they do not fit into Askey's
scheme of hypergeometric orthogonal polynomials. The reason for this is
that the orthogonality measure for the Lommel polynomials is supported on the
set consisting of one over the zeros of a Bessel function, which are not
explicitly known in general. So there is no Rodrigues formula or difference
equation for the Lommel polynomials.

The Bessel function $J_\nu(z)$ of order $\nu$ and argument $z$
is given by the absolutely convergent series expansion
$$
J_\nu(z) = \sum_{k=0}^\infty {{(-1)^k(z/2)^{\nu+2k}}\over{k!\,
\Gamma(\nu+k+1)}}.
\tag{1.1}
$$
The properties of this special function are well understood, see e.g.
the book on Bessel functions by Watson \cite{22}.
A simple recurrence relation for the Bessel functions is,
cf. \cite{22, \S 3.2(1)},
$$
J_{\nu+1}(z) = {{2\nu}\over z}J_\nu(z) - J_{\nu-1}(z).
\tag{1.2}
$$
From iteration of \thetag{1.2} we see that we can express $J_{\nu+m}(z)$
in terms of $J_\nu(z)$ and $J_{\nu-1}(z)$ and the coefficients of
$J_\nu(z)$ and $J_{\nu-1}(z)$ are polynomials in $z^{-1}$. This was
first observed by Lommel in 1871. Explicitly, we have, cf. Watson
\cite{22, \S 9.6},
$$
J_{\nu+m}(z)= h_{m,\nu}({1\over z})J_\nu(z) -
h_{m-1,\nu+1}({1\over z})J_{\nu-1}(z),
\tag{1.3}
$$
where $h_{m,\nu}(z)$ are the Lommel polynomials, which are also known as
associated Lommel polynomials. The Lommel polynomials
satisfy the three-term recurrence relation
$$
h_{m+1,\nu}(z) = 2z(m+\nu) h_{m,\nu}(z)- h_{m-1,\nu}(z),\qquad
h_{-1,\nu}(z)=0, \quad h_{0,\nu}(z)=1.
\tag{1.4}
$$

Favard's theorem, cf. Chihara \cite{7, Ch.~II, thm.~6.4}, implies
that the Lommel polynomials are orthogonal polynomials with respect to
a positive weight function for $\nu>0$. The explicit orthogonality
relations are, cf. Chihara \cite{7, Ch.~VI, \S 6}, Dickinson
\cite{9}, Dickinson, Pollak and Wannier \cite{10}, Ismail
\cite{15}, Schwartz \cite{19},
$$
\sum_{k=1}^\infty {1\over{(j^{\nu-1}_k)^2}}
h_{m,\nu}\Bigl({{\pm 1}\over{j^{\nu-1}_k}}\Bigr)
h_{n,\nu}\Bigl({{\pm 1}\over{j^{\nu-1}_k}}\Bigr)
={{\delta_{n,m}}\over{2(\nu+n)}},
\tag{1.5}
$$
where $j^\nu_k$, $\nu>-1$, are the positive zeros of the Bessel function
$J_\nu(z)$ numbered increasingly, cf. Watson \cite{22, Ch.~15}.

Another relation between the Lommel polynomials and the Bessel function
is given by Hurwitz's asymptotic formula, cf. \cite{22, 9.65(1)},
$$
{{(2z)^{1-\nu-m} h_{m,\nu}(z)}\over{\Gamma(\nu+m)}} \longrightarrow
J_{\nu-1}\Bigl({1\over z}\Bigr),\qquad m\to\infty.
\tag{1.6}
$$

For the Bessel function \thetag{1.1} there exist several $q$-analogues.
The oldest $q$-analogues for the Bessel function were introduced by
Jackson in a series of papers in 1903-1905, see the references in
\cite{15}. For the Jackson $q$-Bessel function Ismail \cite{15}
introduced the associated $q$-Lommel polynomials, which turned out to
satisfy an orthogonality relation similar to \thetag{1.5}, but now
involving the zeros of the Jackson $q$-Bessel function. Ismail used
these $q$-Lommel polynomials
to prove that the zeros of the Jackson $q$-Bessel functions behave
similarly as the zeros of the Bessel function.

A more recent $q$-analogue of the Bessel function has been introduced by
Hahn in a special case and by Exton in full generality, see the
references in Koornwinder and Swarttouw \cite{18}.
The zeros of the Hahn-Exton $q$-Bessel function and several associated
$q$-analogues of the Lommel polynomial have been studied by Koelink and
Swarttouw \cite{17}. The zeros of the Hahn-Exton $q$-Bessel function
behave in a similar fashion as the zeros of the Bessel function.
In that paper \cite{17} a $q$-analogue of the Lommel polynomials was
introduced. However, this $q$-analogue of the Lommel polynomial is
no longer a polynomial, but a Laurent polynomial. One of the
goals of this
paper is to give an explicit orthogonality measure for these orthogonal
Laurent $q$-Lommel polynomials.

The Laurent $q$-Lommel polynomials are defined by, cf. \cite{17,
prop.~4.3 with $R_{m,\nu}(z^{-1};q)=h_{m,\nu}(z;q)$},
$$
h_{m+1,\nu}(x;q) = \bigl( {1\over x}+x(1-q^{\nu+m})\bigr) h_{m,\nu}(x;q)
-h_{m-1,\nu}(x;q),
\tag{1.7}
$$
with initial conditions $h_{-1,\nu}(x;q)=0$, $h_{0,\nu}(x;q)=1$. A
second independent solution of \thetag{1.7} is given by
$h_{m-1,\nu+1}(x;q)$. Note that taking the limit $q\uparrow 1$ in
\thetag{1.7} after replacing $x$ by $2z/(1-q)$ gives \thetag{1.4}.
The Laurent $q$-Lommel polynomials originate from a relation similar
to \thetag{1.3}, see proposition~3.1.

The explicit orthogonality relations for the Laurent $q$-Lommel polynomials
$h_{m,\nu}(x;q)$ defined in \thetag{1.7} is derived in \S 3. The method of
proof is based on the existence of asymptotically well-behaved solutions of
\thetag{1.7} reminiscent of $J_{\nu+m}(x)$, cf. \thetag{1.3}. The method of
Dickinson \cite{9} to prove \thetag{1.5} can then be adapted to our
situation. The orthogonality measure gives rise to a strong moment functional
$\L$, i.e. a functional on the space of Laurent polynomials so that all
moments $\L(x^n)$, $n\in\Z$, exist. From $\L$ we obtain two moment functionals
$\L_{\pm}$, as considered in e.g., \cite{7, Ch.~1}, by putting
$\L_+(x^n)=\L(x^n)$, $n\in\Zp$, and $\L_-(x^n)=-\L(x^{-2-n})$,
$n\in\Zp$. (The $2$ has to do with the fact that all moment functionals
are symmetric.)
It turns out that both $\L_+$ and $\L_-$ are positive definite
moment functionals.

The orthogonal polynomials for $\L_+$ are $q$-analogues of the Lommel
polynomials and the support of orthogonality measure consists of zero and one
over the zeros of a Hahn-Exton $q$-Bessel function. This is worked out in
detail in \S 4, where we use Dickinson's method \cite{9} once more.
In \S 5 we study the orthogonal polynomials for $\L_-$. We give explicit
expressions for these polynomials in terms of Al-Salam--Chihara polynomials,
which can be used to determine the asymptotic behaviour as the degree tends to
infinity. The asymptotic behaviour is expressed in terms of a function
$j_\nu(x;q)$
closely related to the Hahn-Exton $q$-Bessel function. Since we can do this
for the associated polynomials as well, we have the Stieltjes transform of
the orthogonality measure from which the orthogonality follows. Using the
results of \S 5 we can simplify the expression for the strong moment
functional $\L$ using a Wronskian type formula. This is done in \S 6.

For $q=0$, or for $\nu\to\infty$, we see that $U_m\bigl( (x+x^{-1})/2\bigr)$,
where $U_m$ denotes the Chebyshev polynomial of the second kind,
satisfies \thetag{1.7} with the same initial conditions. So we can view the
Laurent $q$-Lommel polynomials $h_{m,\nu}(x;q)$ as a perturbation of the
Chebyshev polynomials. This point of view allows us to obtain estimates for
the Laurent $q$-Lommel polynomials, the Hahn-Exton $q$-Bessel function and
the related function $j_\nu(x;q)$. This is done in \S 7.

Finally, in \S 2 we show that the general theory of orthogonal Laurent
polynomials presents us with an existence theorem for the strong moment
functional $\L$. We also state a result concerning the zeros of the Laurent
$q$-Lommel polynomials.

To end this introduction we briefly recall the notation for basic (or
$q$)-hypergeometric series. We follow the standard notation of Gasper
and Rahman \cite{11, Ch.~1}. We take $0<q<1$ for the rest of the
paper. A $q$-shifted factorial is a product defined by
$$
(a;q)_k = \prod_{i=0}^{k-1}(1-aq^i), \qquad a\in\C, \quad k\in\Zp,
$$
where the empty product equals $1$ by definition. Since $0<q<1$ we can
take $k\to\infty$ to get $\lim_{k\to\infty} (a;q)_k = (a;q)_\infty$.
A basic (or $q$)-hypergeometric series is
$$
\gather
{}_r\vp_s \left( {{a_1,\ldots,a_r}\atop{b_1,\ldots,b_s}};q, z\right)
= {}_r\vp_s (a_1,\ldots,a_r;b_1,\ldots,b_s;q, z) \\
= \sum_{k=0}^\infty {{(a_1;q)_k\ldots (a_r;q)_k}\over{(q;q)_k
(b_1;q)_k\ldots (b_s;q)_k}} \bigl( (-1)^k q^{{1\over
2}k(k-1)}\bigr)^{1+s-r} z^k
\tag{1.8}
\endgather
$$
For generic values of the parameters the radius of convergence of the
series in \thetag{1.8} is $0$, $1$, $\infty$ corresponding to
$r>s+1$, $r=s+1$, $r<s+1$.


\head 2. Orthogonal Laurent polynomials\endhead

In this section we apply some of the theory of orthogonal Laurent
polynomials to the Laurent polynomials $h_{m,\nu}(x;q)$ to obtain the
existence of a strong moment functional $\L$, i.e. a linear functional
on the space of Laurent polynomials for which the moments $\L (x^m)$
exist for all $m\in\Z$, for which the Laurent
$q$-Lommel polynomials are orthogonal. We use
the paper by Hendriksen and van Rossum
\cite{14} as the main reference for this section. The recurrence
relation as in \thetag{2.1} has been generalised to a wider class of
recurrence relations by Ismail and Masson \cite{16} by replacing $x$
in front of the $V_{m-1,\nu}(x)$ by $(x-a_m)$, for which they prove a Favard
type theorem. Specialisation to the case considered here yields the Favard
type theorem contained in Hendriksen and van Rossum \cite{14}.
For further information concerning this section the reader may consult the
introductory paper by Cochran and Cooper \cite{8}.

From the recurrence relation \thetag{1.7} it follows that
$h_{m,\nu}(x;q)$ is an even function for even $m$ and an odd function
for odd $m$.
Consequently, $x^mh_{m,\nu}(x;q)$ is a polynomial in $x^2$, which we denote by
$V_{m,\nu}(x^2)=x^m h_{m,\nu}(x;q)$. For $V_m$ we obtain from \thetag{1.7} the
recurrence relation
$$
V_{m+1,\nu}(x) = \bigl( 1+x(1-q^{\nu+m})\bigr) V_{m,\nu}(x) -
xV_{m-1,\nu}(x)
\tag{2.1}
$$
with initial conditions $V_{-1,\nu}(x)=0$, $V_{0,\nu}(x)=1$,
cf. \cite{14, (2.2)}. The Favard-type
theorem, cf. \cite{14, thm.~1.1},
implies that for the Laurent polynomials $Q_n(x)$, defined by
$$
\align
Q_{2n}(x) &= x^{-n} V_{2n,\nu}(x) = h_{2n,\nu}(\sqrt{x};q) \\
Q_{2n+1}(x) &= x^{-n-1} V_{2n+1,\nu}(x) =
x^{-{1\over 2}}h_{2n+1,\nu}(\sqrt{x};q)
\endalign
$$
there exists a strong moment functional $\L_1$ such that
$\L_1 (Q_nQ_m) =0$ for $n\not= m$.

If we form the lacunary Laurent polynomials, cf. \cite{14,
(1.16)}, we get the Laurent polynomials $P_{2m}(x) = h_{m,\nu}(x;q)$,
$P_{2m+1}(x)=x^{-1}h_{m,\nu}(x;q)$. The lacunary Laurent polynomials are
orthogonal with respect to the strong moment functional $\L$ defined by
$\L (x^{2n})=\L_1(x^n)$, $\L (x^{2n+1})=0$ for $n\in\Z$, cf.
\cite{14, prop.~III}. So the orthogonality relations for the even
lacunary Laurent polynomials gives
$$
\L \bigl( h_{n,\nu}(x;q)h_{m,\nu}(x;q)\bigr)
\cases = 0,&\text{$n\not= m$},\\
   \not= 0,&\text{$n=m$}.\endcases
\tag{2.2}
$$
But we also have the orthogonality for the odd lacunary Laurent
polynomials,
$$
\L \bigl( x^{-2} h_{n,\nu}(x;q)h_{m,\nu}(x;q)\bigr)
\cases = 0,&\text{$n\not= m$},\\
   \not= 0,&\text{$n=m$}.\endcases
\tag{2.3}
$$

The space $\Lambda_n$ of Laurent polynomials of the form
$\sum_{p=-n}^n c_p x^p$ is $(2n+1)$-dimensional, $n\in\Zp$. The Laurent
polynomials $h_{m,\nu}(x;q)$, $m=0,1,\ldots,n$, form a $(n+1)$-dimensional
subspace of $\Lambda_n$. Moreover, they form an orthogonal basis for this
subspace with respect to $\L$. Equation \thetag{2.3} states that this
orthogonal basis can be complemented with $x^{-1}h_{m,\nu}(x;q)$,
$m=0,1,\ldots,n-1$, to give an orthogonal basis of $\Lambda_n$ with
respect to $\L$. $\L(x^{-1}h_{m,\nu}(x;q)h_{n,\nu}(x;q))=0$ is immediate from
$\L (x^{2p+1})=0$.

\demo{Remark~2.1}
For orthonormal polynomials the three-term recurrence relation
can be used to prove that the zeros of the orthonormal polynomials correspond
precisely to the eigenvalues of a truncated Jacobi matrix. A similar approach
can be used here. Define coefficients by
$$
xV_{n,\nu}(x) = \sum_{k=0}^{n+1} c_{n,k} V_{k,\nu}(x),
\tag{2.4}
$$
then the matrix $H_n=(c_{i,j})_{0\leq i,j\leq n-1}$ is Hessenberg matrix,
i.e. $c_{i,j}=0$ for $i+1<j$. Using \thetag{2.4} in \thetag{2.1} gives
recurrence relations for the matrix elements $c_{i,j}$, which can be solved
to give
$$
c_{n,k} = \cases
{\displaystyle{\frac{1}{1-q^{\nu+n}},}} & \text{if $k=n+1$,} \\
{\displaystyle{\frac{(q^\nu;q)_{k-1}q^{\nu+k-1}}{(q^\nu;q)_{n+1}},}}
             & \text{if $0 < k \leq n$,} \\
{\displaystyle{\frac{-1}{(q^\nu;q)_{n+1}},}} & \text{if $k=0$.}
        \endcases
\tag{2.5}
$$
Note that each row sum of $H_n$, except the last, equals zero.

Introduce the vector
$w_n(x) = \bigl( V_{0,\nu}(x), V_{1,\nu}(x), \ldots, V_{n-1,\nu}(x)\bigr)^t$,
then we see from \thetag{2.1} that
$H_n w_n(x) = x w_n(x)$ if $V_{n,\nu}(x)=0$. So a zero $x$ of
$V_{n,\nu}$ implies that $H_n$ has an eigenvector for the eigenvalue $x$.
It is also possible to prove that an eigenvalue $x$ of $H_n$ implies that
$V_{n,\nu}(x)=0$, which can be proved by showing that the characteristic
polynomial of $H_n$ times the normalisation constant $(-1)^n (q^\nu;q)_n$
satisfies \thetag{2.1}. So we conclude that the zeros of $V_{n,\nu}(x)$,
and hence the zeros of the Laurent $q$-Lommel polynomials $h_{n,\nu}(x;q)$,
are completely determined by the spectrum of the Hessenberg matrix $H_n$.
\enddemo


\head 3. Minimal solutions and orthogonality relations\endhead

In this section we give an explicit formulation for the strong moment
functional $\L$ introduced in the previous section. We describe $\L$ in
terms of
contour integrals where the integrands depend on the Hahn-Exton $q$-Bessel
function and on a function closely related to the Hahn-Exton $q$-Bessel
function. These functions give rise to two other solutions of the recurrence
relation \thetag{1.7}, but now with prescribed behaviour for $m\to\infty$.
The proof of
orthogonality of the Laurent $q$-Lommel polynomials for $\L$ uses a method
already introduced by Dickinson \cite{9} to prove the orthogonality
relations \thetag{1.5} for the Lommel polynomials.

Using a generating function argument the following explicit expressions
for the Laurent $q$-Lommel polynomials have been derived in
\cite{17, (4.23)} from the recurrence relation \thetag{1.7}
$$
\align
h_{m,\nu}(x;q) &= \sum_{n=0}^m x^{m-2n}
{{(q^{n+1};q)_\infty (q^\nu;q)_\infty}\over{(q;q)_\infty
(q^{\nu+m-n};q)_\infty}} \, {}_2\vp_1 \left(
{{q^{-n},q^{\nu+m-n}}\atop{q^\nu}};q, q^{n+1}\right) \tag{3.1} \\
&= \sum_{n=0}^m x^{m-2n}
\, {}_2\vp_1 \left({{q^{n-m},q^{n+1}}\atop{q}};q, q^{\nu+m-n}\right).
\tag{3.2}
\endalign
$$

The Hahn-Exton $q$-Bessel function is defined by
$$
J_\nu(x;q) = {{(q^{\nu+1};q)_\infty}\over{(q;q)_\infty}} x^\nu
\, {}_1\varphi_1 \left( {{0}\atop{q^{\nu+1}}};q,qx^2\right)
\tag{3.3}
$$
and then the following $q$-analogue of Hurwitz's formula \thetag{1.6}
holds
$$
\lim_{m\to\infty} x^{-m} h_{m,\nu}(x;q) =
{{(q;q)_\infty}\over{(x^{-2};q)_\infty}} x^{\nu-1} J_{\nu-1}\Bigl({1\over
x};q\Bigr), \qquad \vert x\vert >1.
\tag{3.4}
$$
Relation \thetag{3.4} has been proved formally in Koelink and
Swarttouw \cite{17, (4.24)} from \thetag{3.1}, but it follows from
their proof that it is valid only for $\vert x\vert > 1$.

In order to state the asymptotic behaviour of the Laurent $q$-Lommel
polynomials inside the circle we introduce the function
$$
j_\nu(x;q) = x^\nu (qx^2;q)_\infty \, {}_1\vp_1 (0;qx^2;q,q^{\nu+1}x^2)
= x^\nu (q^{\nu+1}x^2;q)_\infty \, {}_1\vp_1 (0;q^{\nu+1}x^2;q,qx^2),
\tag{3.5}
$$
where we use $(x;q)_\infty \, {}_1
\vp_1 (0;x;q,y) = (y;q)_\infty \,{}_1\vp_1(0;y;q,x)$,
cf. \cite{18, (2.3)}. This function is related to the Hahn-Exton
$q$-Bessel function in the following way
$$
x^{-\nu} j_\nu(x;q) = (q;q)_\infty \Bigl( x^{-\mu} J_\mu(x;q)\Bigr)
\Big\vert_{\mu=\nu+2\ln x/\ln q}.
$$

Now we can use \thetag{3.2} to obtain
$$
\aligned
x^m h_{m,\nu}(x;q) &=
\sum_{n=0}^m x^{2m-2n} \, {}_2\varphi_1 \left( {{q^{n-m},
q^{n+1}}\atop{q}};q,q^{\nu+m-n}\right) \\
& = \sum_{n=0}^m x^{2n} \, {}_2\varphi_1 \left( {{q^{-n},
q^{m-n+1}}\atop{q}};q,q^{\nu+n}\right)
\endaligned
$$
and by dominated convergence we obtain
$$
\aligned
\lim_{m\to\infty} x^m h_{m,\nu}(x;q) & =
\sum_{n=0}^\infty x^{2n} \, {}_2\varphi_1 \left( {{q^{-n},0}
\atop{q}};q,q^{\nu+n}\right) \\
& = \sum_{l=0}^\infty {{q^{\nu l}}\over{(q;q)_l(q;q)_l}}
\sum_{n=l}^\infty (q^{-n};q)_lq^{nl} x^{2n},
\endaligned
$$
where the last equality follows from interchanging the summations, which
is allowed for $\vert x\vert < 1$. The inner sum can be written as
$$
\align
& \sum_{p=0}^\infty (q^{-p-l};q)_l x^{2(p+l)}q^{l(p+l)} =
x^{2l} (-1)^l q^{{1\over 2}l(l-1)}
\sum_{p=0}^\infty (q^{p+l};q^{-1})_l x^{2p} \\
& = x^{2l} (-1)^l q^{{1\over 2}l(l-1)} (q;q)_l
\sum_{p=0}^\infty {{(q^{l+1};q)_p}\over{(q;q)_p}} x^{2p} =
x^{2l} (-1)^l q^{{1\over 2}l(l-1)} {{(q;q)_l}\over{(x^{2};q)_{l+1}}},
\endalign
$$
by the $q$-binomial theorem, cf. \cite{11, (1.3.2)}.
This leads to the result
$$
\lim_{m\to\infty} x^m h_{m,\nu}(x;q) = {1\over{1-x^2}}
\, {}_1\varphi_1 \left( {{0}\atop{qx^2}};q,q^\nu x^2\right)
= {{x^{1-\nu}}\over{(x^2;q)_\infty}} j_{\nu-1}(x;q), \qquad
\vert x\vert < 1.
\tag{3.6}
$$

\proclaim{Proposition~3.1} The functions $J_{\nu+m}(x^{-1};q)$ and
$j_{\nu+m}(x;q)$ satisfy the recurrence relation \thetag{1.7}. Moreover,
$$
\gather
J_{\nu+m}(x^{-1};q) = h_{m,\nu}(x;q) J_\nu (x^{-1};q) -
h_{m-1,\nu+1}(x;q) J_{\nu-1}(x^{-1};q), \\
j_{\nu+m}(x;q) = h_{m,\nu}(x;q) j_\nu (x;q) -
h_{m-1,\nu+1}(x;q) j_{\nu-1}(x;q).
\endgather
$$
\endproclaim

\demo{Proof} Since $h_{m,\nu}(x;q)$ and $h_{m-1,\nu+1}(x;q)$ are
linearly independent solutions of the recurrence relation \thetag{1.7},
the last statement of the proposition implies the first. Also, if
$J_{\nu+m}(x^{-1};q)$ and $j_{\nu+m}(x;q)$ satisfy \thetag{1.7}, then
they must be a linear combination of $h_{m,\nu}(x;q)$ and
$h_{m-1,\nu+1}(x;q)$ from which the second statement follows by
considering the cases $m=0$ and $m=-1$.

The last statement for $J_{\nu+m}$ has already been proved in
\cite{17, (4.12)}, so it remains to consider $j_{\nu+m}$.
The second order
$q$-difference equation for the ${}_1\vp_1$-series, or by
taking a suitable limit in \cite{17, (4.14)}
in combination with \thetag{3.6}, reveals that
$$
j_{\nu+1}(x;q) =
\bigl( {1\over x}+x(1-q^\nu)\bigr) j_\nu(x;q)
-j_{\nu-1}(x;q).
$$
Replacing $\nu$ by $\nu+m$ proves the statement. \qed\enddemo

\demo{Remark} (i) The solutions $J_{\nu+m}(x^{-1};q)$ and
$j_{\nu+m}(x;q)$ of \thetag{1.7} have the following asymptotic behaviour
for $m\to\infty$ valid for $x\in\C$;
$$
\gathered
\lim_{m\to\infty} x^{m+\nu} J_{\nu+m}(x^{-1};q) =
{{(qx^{-2};q)_\infty}\over{(q;q)_\infty}}, \\
\lim_{m\to\infty} x^{-m-\nu} j_{\nu+m}(x;q) = (qx^2;q)_\infty.
\endgathered
\tag{3.7}
$$
Note that $x^{\pm m}$ are solutions of \thetag{1.7} for $m\to\infty$
(or for $q=0)$. So the solutions $J_{\nu+m}(x^{-1};q)$
and $j_{\nu+m}(x;q)$ behave as $x^{\mp m}$ up to a factor independent of $m$
as $m\to\infty$.

\noindent
(ii) The functions $J_{\nu+m}(x^{-1};q)$ and $j_{\nu+m}(x;q)$ are
related to a minimal solution $X_m(x)$ of \thetag{2.1}, i.e. $X_m(x)$ is
a solution such that $\lim_{m\to\infty} X_m(x)/V_{m,\nu}(x)=0$,
where $V_{m,\nu}(x)$
is the polynomial solution of \thetag{2.1}. Using the limit transitions
\thetag{3.4} and \thetag{3.6} and the relations in proposition~3.1 we
obtain
$$
X_m(x) = \cases j_\nu(\sqrt{x};q) V_{m,\nu}(x) -
x^{1\over 2}j_{\nu-1}(\sqrt{x};q) V_{m-1,\nu+1}(x) &\text{}\\
\qquad\qquad\qquad\qquad\qquad\qquad\qquad= x^{{1\over 2}m}
j_{\nu+m}(\sqrt{x};q),&\text{$\vert x\vert <1$}, \\
J_\nu(1/\sqrt{x};q) V_{m,\nu}(x) -
x^{1\over 2}J_{\nu-1}(1/\sqrt{x};q) V_{m-1,\nu+1}(x) &\text{}\\
\qquad\qquad\qquad\qquad\qquad\qquad\qquad = x^{{1\over 2}m}
J_{\nu+m}(1/\sqrt{x};q),&\text{$\vert x\vert >1$}.
\endcases
$$
\enddemo

With the functions $J_\nu(x;q)$ and
$j_\nu(x;q)$ and their relation with the Laurent $q$-Lommel polynomials
described in proposition~3.1 at hand, we can give an explicit expression
for the strong moment functional $\L$. The proof we give is an adaption to
the Laurent case of Dickinson's proof of the orthogonality \thetag{1.5}
for the Lommel polynomials \cite{9}.

First we investigate the quotient of two Hahn-Exton $q$-Bessel
functions.

\proclaim{Lemma~3.2}
For $\nu>0$ the following expansion holds around $0$ for $n\in\Zp$
$$
{{J_{\nu+n}(x;q)}\over{J_{\nu-1}(x;q)}} =
{{x^{n+1}}\over{(q^\nu;q)_{n+1}}}\sum_{k=0}^\infty c_k x^{2k}
$$
with the coefficients $c_k$ recursively defined by
$c_0=1$ and
$$
c_k = {{(-1)^k q^{{1\over 2}k(k+1)}}\over{(q^{\nu+n+1};q)_k(q;q)_k}}
- \sum_{p=0}^{k-1} c_p
{{(-1)^{k-p} q^{{1\over 2}(k-p)(k-p+1)}}\over{(q^\nu;q)_k(q;q)_{k-p}}}.
\tag{3.8}
$$
\endproclaim

\demo{Proof} From \thetag{3.3} we immediately get
$$
{{J_{\nu+n}(x;q)}\over{J_{\nu-1}(x;q)}} =
{{x^{n+1}}\over{(q^\nu;q)_{n+1}}}
{{ {}_1\varphi_1(0;q^{\nu+n+1};q,qx^2)}\over{
{}_1\varphi_1(0;q^\nu;q,qx^2)}}
$$
so we have to solve for the coefficients $c_k$ by comparing powers of $x$ on
both sides of
$$
\sum_{k=0}^\infty c_k x^{2k} \sum_{p=0}^\infty {{(-1)^p q^{{1\over
2}p(p+1)} x^{2p}}\over{(q^\nu;q)_p(q;q)_p}} =
\sum_{m=0}^\infty {{(-1)^m q^{{1\over
2}m(m+1)} x^{2m}}\over{(q^{\nu+n+1};q)_m(q;q)_m}}
$$
from which the recurrence relation \thetag{3.8}
for the coefficients $c_k$ is obtained.

A rude estimate gives
$$
\Bigl\vert {{(-1)^k q^{{1\over 2}k(k+1)}}\over{(q^{\nu+n+1};q)_r(q;q)_r}}
\Bigr\vert \leq A= {1\over{(q^\nu;q)_\infty (q;q)_\infty}}
$$
for $\nu>0$. The same estimate applies to the factor in front of $c_n$
on the right hand side of \thetag{3.8}, so that we obtain
$$
\vert c_k\vert \leq A + \sum_{n=0}^{k-1} A \vert c_n\vert.
$$
A discrete version of Gronwall's inequality, cf.
e.g. \cite{21, p.~440}
$$
a_k\leq A + \sum_{n=0}^{k-1} d_n a_n, \ \hbox{\rm and}\
A,a_n,d_n\geq 0 \Longrightarrow a_n \leq A\exp \Bigl( \sum_{n=0}^{k-1}
d_n\Bigr),
\tag{3.9}
$$
yields $\vert c_k\vert \leq Ae^{kA}$ so that the series on the right hand
side of the statement of the lemma
is absolutely convergent for $\vert x\vert < e^{-A/2}$.
\qed\enddemo

Choose $0<R< j^{\nu-1}_1$, where $j_1^{\nu-1}$ denotes the smallest
positive zero of $J_{\nu-1}(x;q)$, $\nu>0$, cf. \cite{17, \S 3}.
Using lemma~3.2  we obtain for $\nu>0$, $m\in\Z$ and $n\in\Zp$
$$
{1\over{2\pi i}} \oint_{\vert z\vert ={1\over R}} z^m
{{J_{\nu+n}(z^{-1};q)}\over{J_{\nu-1}(z^{-1};q)}}\, dz =
\cases 0, &\text{$m-n$ odd or $m<n$,} \\
   (q^\nu;q)_{n+1}^{-1}, &\text{$m=n$.} \endcases
\tag{3.10}
$$
Note that the coefficients $c_k$ of lemma~3.2 for $n=0$
are in fact the moments of the linear
functional $\L_+$ defined by
$$
\L_+(x^m) =
{1\over{2\pi i}} \oint_{\vert z\vert ={1\over R}} z^m
{{J_\nu(z^{-1};q)}\over{J_{\nu-1}(z^{-1};q)}}\, dz =
\cases 0, &\text{$m\in\Zp$ odd,} \\
       c_{m/2}, &\text{$m\in\Zp$ even.} \endcases
\tag{3.11}
$$
We return to this moment functional in section 4 and we calculate
the corresponding orthogonal polynomials, which turn out to be
$q$-analogues of the Lommel polynomials.

The following lemma is the analogue of lemma~3.2 for the functions
$j_\nu(x;q)$ instead of the Hahn-Exton $q$-Bessel function.

\proclaim{Lemma~3.3}
For $\nu\in\R$ the following expansion holds around $0$ for $n\in\Zp$
$$
{{j_{\nu+n}(x;q)}\over{j_{\nu-1}(x;q)}} =
x^{n+1}\sum_{k=0}^\infty d_k x^{2k}
$$
with the coefficients $d_k$ recursively defined by
$d_0=1$ and
$$
d_k = {}_2\varphi_1(q^{-k},0;q;q,q^{\nu+n+1+k}) -
 \sum_{p=0}^{k-1} d_p\ {}_2\varphi_1(q^{p-k},0;q;q,q^{\nu+k-p})
\tag{3.12}
$$
\endproclaim

\demo{Proof} The proof is completely analogous to the proof of lemma~3.2
and we only give the differences. Here use the expansion
$$
\aligned
&{1\over{(1-x^2)}} \, {}_1\varphi_1 (0;qx^2;q,q^{\nu+1}x^2) =
\sum_{k=0}^\infty {{(-1)^k q^{{1\over
2}k(k-1)}}\over{(q;q)_k(x^2;q)_{k+1}}} q^{(\nu+1)k} x^{2k} \\
&\qquad = \sum_{k=0}^\infty \sum_{l=0}^\infty
{{(-1)^k q^{{1\over 2}k(k-1)}}\over{(q;q)_k}}
q^{(\nu+1)k} x^{2k} {{(q^{k+1};q)_l}\over{(q;q)_l}} x^{2l} \\
&\qquad = \sum_{p=0}^\infty x^{2p} \, {}_2\varphi_1
(q^{-p},0;q;q,q^{\nu+1+p})
\endaligned
$$
by the $q$-binomial theorem valid for $\vert x\vert < 1$ and
rearranging the absolutely convergent sum using $l=p-k$.
From this we obtain the
recurrence relation \thetag{3.12}.
The general estimate
$$
\vert {}_2\varphi_1(q^{-p},0;q;q,xq^p)\vert \leq {{(-q,-\vert
x\vert;q)_\infty}\over{(q;q)_\infty}}
$$
and Gronwall's inequality \thetag{3.9} prove that sum is absolutely
convergent around $0$.
\qed\enddemo

Choose $r>0$ so that $j_{\nu-1}(x;q)$ has no non-zero
zeros in the ball with
radius $r$ and the origin as centre, which is possible since
${}_1\varphi_1(0;qx^2;q,q^\nu x^2)$ equals $1$ at $x=0$ and defines an
analytic function for $\vert x\vert <q^{-1/2}$.
Using lemma~3.3 we obtain for $m\in\Z$ and $n\in\Zp$
$$
{1\over{2\pi i}} \oint_{\vert z\vert =r} z^m
{{j_{\nu+n}(z;q)}\over{j_{\nu-1}(z;q)}}\, dz =
\cases 0, &\text{$m-n$ odd or $m>-n-2$,} \\
   1, &\text{$m=-n-2$.} \endcases
\tag{3.13}
$$

The coefficients $d_k$ of lemma~3.3 for $n=0$ can be interpreted as the
moments of the moment functional $\L_-$ defined by
$$
\L_-(x^m) = {1\over{2\pi i}} \oint_{\vert z\vert = {1\over r}} z^m
{{j_\nu(z^{-1};q)}\over{j_{\nu-1}(z^{-1};q)}} \, dz =
\cases 0, &\text{$m\in\Zp$ odd,} \\
       d_{m/2}, &\text{$m\in\Zp$ even.} \endcases
\tag{3.14}
$$
In \S 5 we consider the orthogonal polynomials for $\L_-$ from which
some properties for $j_\nu(x;q)$ can be derived.

Define the following strong moment functional $\L$ for $\nu>0$
on the space of Laurent polynomials by
$$
\L(p) = {1\over{2\pi i}} \oint_{\vert z\vert = {1\over R}}
p(z) {{J_\nu(z^{-1};q)}\over{J_{\nu-1}(z^{-1};q)}}\, dz -
{1\over {2\pi i}} \oint_{\vert z\vert = r}
p(z) {{j_\nu(z;q)}\over{j_{\nu-1}(z;q)}}\, dz
\tag{3.15}
$$
for any Laurent polynomial $p(z) = \sum_{p=n}^m c_p z^p$, $n\leq m$,
$n,m\in\Z$. Note that $\L$ is independent of the choice of $R$, respectively
$r$, as long as $J_{\nu-1}(x;q)$, respectively $j_{\nu-1}(x;q)$, has no
non-zero zeros in the ball with radius $R$, respectively $r$.
All moments of $\L$, both positive and negative, are well-defined
due to lemmas~3.2 and 3.3.

The moments of strong moment functional $\L$ and the moments of the moment
functionals $\L_\pm$ defined in \thetag{3.11} and \thetag{3.14} are related
by $\L_+(x^n)=\L(x^n)$, $n\in\Zp$, and by
$\L_-(x^n)=-\L(x^{-2-n})$, $n\in\Zp$.

\proclaim{Theorem~3.4}
Let $\nu>0$.
The Laurent $q$-Lommel polynomials $h_{n,\nu}(x;q)$ defined by \thetag{1.7}
are orthogonal Laurent polynomials with respect to the strong moment
functional $\L$, cf. \thetag{3.15}.
Moreover, also the Laurent polynomials $x^{-1}h_{n,\nu}(x;q)$ are orthogonal
with respect to $\L$. Explicitly,
$$
\L (h_{n,\nu}(x;q) h_{m,\nu}(x;q)) =
{{\delta_{n,m}}\over{1-q^{\nu+m}}}, \qquad
\L (x^{-1}h_{n,\nu}(x;q) x^{-1}h_{m,\nu}(x;q)) =
- \delta_{n,m}.
$$
\endproclaim

\demo{Remark} (i) This result corresponds nicely with the fact that the
Laurent $q$-Lommel polynomials correspond to a sequence of lacunary
orthogonal Laurent polyomials, cf. \thetag{2.2}, \thetag{2.3}.

\noindent
(ii) Since $\L(x^{-2}) = -1$ we see that $\L$ is not a
positive definite strong moment functional.
\enddemo

\demo{Proof} The asymptotically
well-behaved solutions $J_{\nu+n}(x^{-1};q)$ and $j_{\nu+n}(x;q)$
of the recurrence relation \thetag{1.7} are expressible in
terms of the Laurent polynomials $h_{n,\nu}(x;q)$ and the associated
Laurent polynomials $h_{n-1,\nu+1}(x;q)$, cf. proposition~3.1.
From this we obtain for any $m\in\Z$ the expressions
$$
\gather
x^m {{J_{\nu+n}(x^{-1};q)}\over{J_{\nu-1}(x^{-1};q)}} =
x^m {{J_\nu(x^{-1};q)}\over{J_{\nu-1}(x^{-1};q)}} h_{n,\nu}(x;q) -
x^m h_{n-1,\nu+1}(x;q), \tag{3.16} \\
x^m {{j_{\nu+n}(x;q)}\over{j_{\nu-1}(x;q)}} =
x^m {{j_\nu(x;q)}\over{j_{\nu-1}(x;q)}} h_{n,\nu}(x;q) -
x^m h_{n-1,\nu+1}(x;q). \tag{3.17}
\endgather
$$

Since we obviously have
$$
{1\over{2\pi i}} \oint_{\vert z\vert={1\over R}} z^m h_{n-1,\nu+1}(z;q)\,
dz = {1\over{2\pi i}} \oint_{\vert z\vert=r} z^m h_{n-1,\nu+1}(z;q)\, dz,
$$
we get from the combination of \thetag{3.16}, \thetag{3.17},
\thetag{3.10} and \thetag{3.13} the relations
$$
\gather
\L(x^m h_{n,\nu}(x;q)) = \cases 0, &\text{$-n\leq m <n$,}\\
 (q^\nu;q)_{n+1}^{-1}, &\text{$m=n$,}\endcases \\
\L(x^m x^{-1}h_{n,\nu}(x;q)) = \cases 0, &\text{$-n<m\leq n$,}\\
 -1, &\text{$m=-n-1$,}\endcases
\endgather
$$
This proves the orthogonality.

It remains to calculate the norm. From \thetag{3.1}, \thetag{3.2}
we see that the coefficient of $x^n$ in
$h_{n,\nu}(x;q)$ equals $(q^\nu;q)_n$ and that the coefficient of
$x^{-n-1}$ in $x^{-1}h_{n,\nu}(x;q)$ equals $1$.
\qed\enddemo


\head 4. Orthogonal $q$-Lommel polynomials associated with the
positive moments\endhead

In this section we consider the orthogonal polynomials for the moment
functional $\L_+$, cf. \thetag{3.11},
which corresponds to the positive moments of
the strong moment functional $\L$. These polynomials are
$q$-analogues of the Lommel polynomials $h_{n,\nu}(z)$, cf.
\thetag{1.4}.

We consider the following three-term recurrence relation,
$$
p_{n+1}(x) = x(1-q^{\nu +n}) p_n(x) - \lambda_n p_{n-1}(x),
\qquad \lambda_{2n}= q^n,\quad \lambda_{2n+1}=q^{\nu+3n+1}
\tag{4.1}
$$
with initial conditions $p_{-1}(x)=0$, $p_0(x) = 1$. Note that we can write
the recurrence coefficient $\lambda_n$ in closed form as
$q^{(\nu +n)(\lfloor (n+1)/2\rfloor - \lfloor n/2\rfloor) +
\lfloor n/2\rfloor}$,  where $\lfloor a\rfloor$ denotes the greatest integer
less than or equal to $a\in\R$.
So the recurrence relation \thetag{4.1}
depends on whether $n$ is odd or even. Favard's theorem implies that
these polynomials are orthogonal with respect to a positive definite moment
functional for $\nu>0$. Taking $q\uparrow 1$ in
\thetag{4.1} after replacing $x$ by $2z/(1-q)$
yields the three-term recurrence relation
\thetag{1.4} for the Lommel polynomials, so that we have
$q$-analogues of the Lommel polynomials.
The recurrence relation \thetag{4.1} has been guessed
using the explicit form for the positive moments of $\L$,
i.e. the moments of $\L_+$,
obtainable from lemma~3.2 and calculating the first
few terms of the recurrence relation \thetag{4.1}
using Mathematica.

The monic orthogonal polynomials satisfy a recurrence relation of the type
$$
r_{n+1}(x) = xr_n(x) - \mu_n r_{n-1}(x), \qquad r_{-1}(x)=0,\ r_0(x)=1,
$$
with $\mu_n>0$ for all $n\in\N$ and $\sum_{n=1}^\infty \mu_n<\infty$.
This type of orthogonal polynomials has been studied by Dickinson, Pollak
and Wannier \cite{10}, by Goldberg \cite{13}, who corrected some of
the results of \cite{10}, and, from the point of view of continued
fractions, by Schwartz \cite{19}. See also Chihara \cite{7, Ch.~IV,
thm.~3.5}. The support of the corresponding orthogonality measure, which is
uniquely determined, is a purely discrete denumerable bounded set with only
one accumulation point at zero. This result can also be obtained by remarking
that the Jacobi matrix $J$ for the corresponding orthonormal polynomials
defines a self-adjoint operator
$J\colon\ell^2(\Zp)\to\ell^2(\Zp)$, which is an operator of trace class.
Since the spectral measure of $J$ is the orthogonality measure for the
orthogonal polynomials $r_n$, the result follows from standard facts on the
spectral measure of a self-adjoint trace class operator. Moreover, for the
orthogonal polynomials in this class we have the asymptotic behaviour of the
form $\lim_{n\to\infty} x^{-n}r_n(x)=f(x)$ for an analytic function $f$ in
$\C\backslash \{ 0\}$, cf. \cite{10}, \cite{13}, \cite{19}.

By $p_n^{(1)}$ we denote the
associated orthogonal polynomials, i.e. the polynomials satisfying
$$
p_n^{(1)}(x) = x(1-q^{\nu +n}) p_{n-1}^{(1)}(x) - \lambda_n
p_{n-2}^{(1)}(x),
\qquad \lambda_{2n}= q^n,\quad \lambda_{2n+1}=q^{\nu+3n+1}
\tag{4.2}
$$
with initial conditions $p^{(1)}_{-1}(x)=0$, $p_0^{(1)}(x) = 1$.

The following proposition is a $q$-analogue of the identity \thetag{1.3}
relating the Bessel functions and Lommel polynomials.

\proclaim{Proposition~4.1} For $n\in\Zp$
the polynomials defined by
\thetag{4.1} and
\thetag{4.2} satisfy
$$
p_n\bigl({1\over x}\bigr) J_\nu(x;q) -
p^{(1)}_{n-1}\bigl({1\over x}\bigr) J_{\nu -1}(x;q)
= q^{{1\over 2}\lfloor (n+1)/2\rfloor(n+\nu)} J_{\nu+n}(xq^{{1\over
2}\lfloor (n+1)/2\rfloor};q)
$$
where $J_\nu(x;q)$ denotes the Hahn-Exton $q$-Bessel function
\thetag{3.3}.
\endproclaim

\demo{Proof} The left hand side is a solution of the three-term recurrence
relation
$$
a_{n+1} = {{1-q^{\nu+n}}\over x} a_n - \lambda_n a_{n-1}.
\tag{4.3}
$$
The right hand side satisfies the same recurrence relation \thetag{4.3}. To
see this we use for even $n$ the relation
$$
{{1-q^\nu}\over x} J_\nu(x;q) - J_{\nu-1}(x;q) = q^{{1\over 2}(\nu+1)}
J_{\nu+1}(xq^{1\over 2};q),
\tag{4.4}
$$
and for odd $n$ we use the relation
$$
{{1-q^\nu}\over x} J_\nu(x;q) - q^{{1\over 2}(\nu-1)}
J_{\nu-1}(xq^{-{1\over 2}};q) = J_{\nu+1}(x;q).
\tag{4.5}
$$
These identities can be checked straightforwardly by comparing the
coefficients of the powers of $x$ on both sides of \thetag{4.4} and
\thetag{4.5}.

Since $p_n(x^{-1})$ and $p_{n-1}^{(1)}(x^{-1})$ are linearly independent
solutions of \thetag{4.3} we obtain the proposition after checking
the equality for $n=0$, which is trivial, and for $n=1$, which is
\thetag{4.4}.
\qed\enddemo

The polynomials defined by \thetag{4.1} turn out to be
the orthogonal polynomials with respect to the moment functional
$\L_+$. For more information concerning the zeros of the Hahn-Exton
$q$-Bessel function, which play a role in the following theorem, the reader
is referred to \cite{17, \S 3}.

\proclaim{Theorem~4.2} We have
the following orthogonality relations for $\nu>0$ for the
polynomials defined by \thetag{4.1};
$$
\gather
\sum_{k=1}^\infty
p_n\Bigl({{\pm 1}\over{j_k^{\nu-1}}}\Bigr)
p_m\Bigl({{\pm 1}\over{j_k^{\nu-1}}}\Bigr)
{{-J_\nu(j_k^{\nu-1};q)}\over{
(j_k^{\nu-1})^2 J^\prime_{\nu-1}(j_k^{\nu-1};q)}} \\
+  p_n(0)p_m(0) =
\delta_{n,m} {{q^{(n+\nu)\lfloor (n+1)/2\rfloor}}\over{1-q^{n+\nu}}}.
\endgather
$$
Here $j_k^{\nu-1}$ are the positive simple zeros of
the Hahn-Exton $q$-Bessel function $J_{\nu-1}(x;q)$
numbered increasingly. All weights are positive.
\endproclaim

\demo{Proof} We start by establishing a complex orthogonality similarly to the
previous section following Dickinson's method \cite{9}.
For this we need the expansion
$$
q^{{1\over 2}\lfloor (n+1)/2\rfloor(n+\nu)} {{J_{\nu+n}(xq^{{1\over
2}\lfloor (n+1)/2\rfloor};q)}\over{J_{\nu-1}(x;q)}} =
q^{\lfloor (n+1)/2\rfloor(n+\nu)} {{x^{n+1}}\over{(q^\nu;q)_{n+1}}}
\sum_{k=0}^\infty c_k x^{2k},
\tag{4.6}
$$
which is absolutely convergent for small $x$. Moreover, $c_0=1$. This is
proved similarly as in lemmas~3.2 and 3.3.

Let $R>0$ be smaller than the smallest positive zero $j_1^{\nu-1}$ of
$J_{\nu-1}(x;q)$, then we obtain from proposition~4.1
and \thetag{4.6} for $0\leq m \leq n$
$$
\aligned
\oint_{\vert z\vert ={1\over R}} z^m p_n(z)
{{J_\nu(z^{-1};q)}\over{J_{\nu-1}(z^{-1};q)}}dz & =
\oint_{\vert z\vert = {1\over R}} z^m
q^{{1\over 2}\lfloor (n+1)/2\rfloor(n+\nu)}
{{J_{\nu+n}(z^{-1}q^{{1\over 2}\lfloor (n+1)/2\rfloor};q)}
\over{J_{\nu-1}(z^{-1};q)}}dz \\
& = \cases 0, &\text{$0\leq m<n$,} \\
2\pi i
q^{(\nu+n)\lfloor (n+1)/2\rfloor} (q^\nu;q)_{n+1}^{-1}, &\text{$m=n$,}
\endcases
\endaligned
$$
since $\oint_{\vert z\vert ={1\over R}} z^m
p_{n-1}^{(1)}(z) dz=0$. The leading coefficient of $p_n$ is
$(q^\nu;q)_n$, as can be read off from \thetag{4.1}, and
so we get the complex orthogonality relations
$$
\L_+(p_np_m) =
{1\over{2\pi i}}\oint_{\vert z\vert ={1\over R}} p_m(z) p_n(z)
{{J_\nu(z^{-1};q)}\over{J_{\nu-1}(z^{-1};q)}}dz  =
\delta_{n,m} {{q^{(n+\nu)\lfloor (n+1)/2\rfloor}}\over{(1-q^{\nu+n})}}.
\tag{4.7}
$$

The considerations given at the beginning of this section show that
we can rewrite \thetag{4.7} as a
sum over the zeros of the Hahn-Exton $q$-Bessel function
$J_{\nu-1}(z;q)$ and possibly zero. The residues at the pole
$(j_k^{\nu-1})^{-1}$ of the left hand side of \thetag{4.7}
equals
$$
p_n({1\over{j_k^{\nu-1}}})
p_m({1\over{j_k^{\nu-1}}})
{{-J_\nu(j_k^{\nu-1};q)}\over{
(j_k^{\nu-1})^2 J^\prime_{\nu-1}(j_k^{\nu-1};q)}}
$$
To see this we note that $J^\prime_{\nu-1}(j_k^{\nu-1};q)\not =0$ since
the zeros of $J_{\nu-1}(x;q)$ are simple, cf. \cite{17, lemma~3.3},
and that $J_\nu(j_k^{\nu-1};q)\not= 0$ by the interlacing
property of the zeros of the Hahn-Exton $q$-Bessel function \cite{17,
thm.~3.7}. The positivity of the corresponding mass follows from the fact
that $J_\nu(j_k^{\nu-1};q)$ and
$J^\prime_{\nu-1}(j_k^{\nu-1};q)$ have opposite signs, which follows
from the Fourier-Bessel orthogonality
relations for the Hahn-Exton $q$-Bessel function \cite{17, prop.~3.6},
or from the fact that the zeros of the Hahn-Exton $q$-Bessel functions
$J_\nu(x;q)$ and $J_{\nu+1}(x;q)$ are interlaced as described in
\cite{17, thm.~3.7}. The mass at
$-(j_k^{\nu-1})^{-1}$ yields the same weight.

The set of mass points $(j_k^{\nu-1})^{-1}$, $k\in\N$, has zero as the
only point of accumulation, so that zero may occur as a mass point as
well. This happens if $\sum_{k=0}^\infty \vert \tilde p_k(0)\vert^2
< \infty$, where $\tilde p_n$ are the corresponding orthonormal
polynomials, cf. e.g. \cite{5, thm.~2.8}.
Now the orthonormal polynomials $\tilde p_n$ are given by
$$
\tilde p_n(x) = \left( {{1-q^{n+\nu}}\over{1-q^\nu}} \right)^{1/2}
q^{-{1\over 2}\lfloor (n+1)/2\rfloor (n+\nu)} p_n(x).
\tag{4.8}
$$
Moreover, ${\Cal M}(\tilde p_n\tilde p_m) = \delta_{n,m}$, where ${\Cal
M}$ is the moment functional given by
$$
{\Cal M}(p) = {{1-q^\nu}\over{2\pi i}}
\oint_{\vert z\vert ={1\over R}} p(z)
{{J_\nu(z^{-1};q)}\over{J_{\nu-1}(z^{-1};q)}}\, dz
= (1-q^\nu) \L_+(p).
$$

From \thetag{4.1} with $x=0$ we see that $p_{2n+1}(0)=0$ and that $p_{2n}(0)$
satisfies a simple two-term recurrence relation from which we get
$p_{2n}(0)=(-1)^n q^{n(\nu+1) + {3\over 2}n(n-1)}$. Combining this with
\thetag{4.8} shows that $\tilde p_{2n+1}(0)=0$ and
$$
\tilde p_{2n}(0) = (-1)^n \left( {{1-q^{\nu+2n}}\over{1-q^\nu}}\right)^{1/2}
q^{{1\over 2}n\nu + {1\over 2}n(n-1)}.
$$
Hence,
$$
\rho = \sum_{k=0}^\infty \vert \tilde p_k(0)\vert^2 = {1\over{1-q^\nu}}
\sum_{n=0}^\infty (1-q^{\nu+2n}) q^{n\nu + n(n-1)}
$$
and this sum is an absolutely convergent
telescoping series so that $\rho = (1-q^\nu)^{-1}$.
Consequently, ${\Cal M}$ has a mass point at zero with weight
$\rho^{-1}$ and $\L_+$ has a mass point at zero with weight $1$.
\qed\enddemo

From the explicit orthogonality relations of theorem~4.2 we see that the
orthogonality measure for $p_n(x)$ is supported in
$[ -1/j_1^{\nu-1},1/j_1^{\nu-1}]$. On the other hand, from the explicit values
of the recurrence coefficients for the orthonormal polynomials $\tilde p_n$,
which are easily obtained from \thetag{4.1} and \thetag{4.8}, and the bound
on the spectrum from \cite{21, (1.3) with $n\to\infty$}, which is
Gershgorin's theorem for the Jacobi matrix, we see that the
orthogonality measure is supported in $[-N,N]$ with $N\leq 2/(1-q^\nu)$.
So we obtain the following corollary after
shifting $\nu$ by $1$.

\proclaim{Corollary~4.3} For $\nu>-1$ the first positive zero $j_1^\nu$
of $J_\nu(x;q)$ satisfies $j_1^\nu \geq (1-q^{\nu+1})/2$.
\endproclaim


\head 5. Orthogonal polynomials associated with the
negative moments\endhead

In this section we consider the orthogonal polynomials for the moment
functional $\L_-$ related to the negative moments of the strong moment
functional $\L$ introduced in \thetag{3.14}.
In subsection~5.1 we introduce the three-term recurrence relation for the
polynomials which we study. The three-term recurrence relation has been
obtained by calculating the first few recurrence coefficients using lemma~3.3
with $n=0$ using Mathematica and then guessing the general result. In
subsection~5.1 we give explicit expressions for these orthogonal polynomials
and the associated orthogonal polynomials in terms of
Al-Salam--Chihara polynomials. From the explicit expressions we can determine
the asymptotic behaviour of the (associated) polynomials as the degree tends
to infinity in terms of the function $j_\nu(x;q)$. In particular, we obtain
the Stieltjes transform of the orthogonality measure. In subsection~5.2
we use the Stieltjes transform to obtain information on the zeros of
$j_\nu(x;q)$ in a similar way as in Ismail \cite{15, \S 4}, see also
\cite{2, \S 4}, and to give explicit
orthogonality relations. In subsection~5.3 we give a different derivation
of some of these results in the special case $\nu=1/2$, which turns out to be
related to known orthogonal polynomials \cite{2}, \cite{20}.
Comparison of these two approaches yields a summation formula for a
one-parameter terminating ${}_3\vp_2$-series.
A special case of this summation formula
is the evaluation of the continuous
$q$-Hermite polynomials at a special point outside the spectrum.

\subhead 5.1. Explicit expressions for orthogonal polynomials\endsubhead
We investigate the monic orthogonal polynomials satisfying the
three-term recurrence relation
$$
P_{n+1}(x) = xP_n(x) - \lambda_n P_{n-1}(x), \qquad
\lambda_{2n} = q^n, \quad  \lambda_{2n+1} = q^{n+\nu},
\tag{5.1}
$$
with initial conditions $P_{-1}(x)=0$, $P_0(x)=1$. By Favard's theorem the
polynomials $P_n$ are orthogonal with respect to a positive definite moment
functional for $\nu\in\R$. Moreover, the polynomials $P_n$ fit into the same
class of Dickinson, Pollak and Wannier \cite{10}, Goldberg \cite{13}
and Schwartz \cite{19} described at the beginning of the previous section.

The polynomials $P_n$ are
even functions of $x$ for even $n$, and odd functions of $x$ for
odd $n$. Introduce
$$
P_{2n}(x) = R_n(x^2), \quad P_{2n+1}(x) = x S_n(x^2),
$$
so that the monic polynomials $R_n$ and $S_n$ satisfy
the three-term recurrence
relations (see, e.g., Chihara \cite{7, p.~45})
$$
\align
R_{n+1}(x) &= (x-\lambda_{2n}-\lambda_{2n+1})R_n(x) -
     \lambda_{2n-1}\lambda_{2n} R_{n-1}(x), \\
S_{n+1}(x) &= (x-\lambda_{2n+1}-\lambda_{2n+2})S_n(x) -
     \lambda_{2n}\lambda_{2n+1} S_{n-1}(x),
\endalign
$$
with initial conditions $R_0(x)=1$, $R_1(x)=x-q^\nu$ and $S_{-1}(x)=0$,
$S_0(x)=1$.
A simple computation from \thetag{5.1} gives the recurrence coefficients
for the polynomials $R_n$
$$
\lambda_{2n}+\lambda_{2n+1} = \cases (1+q^\nu) q^n & \text{if $n>0$}, \\
                                      q^\nu & \text{if $n=0$.}
                               \endcases , \quad
\lambda_{2n}\lambda_{2n-1} = q^{2n-1+\nu}, \qquad n \geq 0.
$$
For the recurrence coefficients of $S_n$ we find similarly
$$
\lambda_{2n+1}+\lambda_{2n+2} = (q+q^\nu) q^n, \quad
\lambda_{2n}\lambda_{2n+1} = q^{2n+\nu}, \qquad n \geq 0.
$$
The recurrence coefficients of $R_n$ and $S_n$ decrease exponentially.

Consider the monic polynomials $u_n(x;a,b;q)$ satisfying the recurrence
relation
$$
u_{n+1}(x;a,b;q) = (x-aq^n)u_n(x;a,b;q) - b^2 q^{2n-2}u_{n-1}(x;a,b;q),
\tag{5.2}
$$
$u_{-1}(x)=0$, $u_0(x)=1$, which are studied in \cite{20}, then
$S_n(x) = u_n(x;q+q^\nu,q^{(\nu+2)/2};q)$. For $R_n$ we have to be a little
bit more careful, since for $n=0$ one of the recurrence
coefficients behaves differently. However, $R_n$ is still a solution of the
recurrence relation \thetag{5.2} with $a=1+q^\nu$ and $b^2=q^{\nu+1}$,
but it satisfies the different
initial condition $R_1(x)=x-q^\nu=u_1(x)+1$. Such polynomials are known as
co-recursive polynomials \cite{6} and can be expressed as
$$
R_n(x) = u_n(x;1+q^\nu,q^{(\nu+1)/2};q) +
u_{n-1}^{(1)}(x;1+q^\nu,q^{(\nu+1)/2};q).
$$
The associated polynomials corresponding to the recurrence
relation \thetag{5.2}
are given by $u_n^{(1)}(x;a,b;q)=u_n(x;aq,bq;q) = q^n u_n(x/q;a,b;q)$,
so that
$$
R_n(x) = u_n(x;1+q^\nu,q^{(\nu+1)/2};q) +
q^{n-1} u_{n-1}(x/q;1+q^\nu,q^{(\nu+1)/2};q).
$$

An explicit expression of the polynomials $u_n(x;a,b;q)$ in terms of
Al-Salam--Chihara polynomials is given by Van Assche \cite{20, thm.~2};
$$
u_n(x;a,b;q) = \sum_{k=0}^n {{x^{n-k}q^{{1\over 2}k(k-1)}}\over{(q;q)_k}}
P_k(-a;q;-aq^{n-k+1},b^2q^{2(n-k)+1},b^2/q).
$$
Here $P_n(x;q;a,b,c)$ are Al-Salam--Chihara polynomials, cf. \cite{1,
(6.1)}, which satisfy the recurrence relation
$$
P_{n+1}(x;q;a,b,c) = (x-aq^n) P_n(x;q;a,b,c) - (c-bq^{n-1})(1-q^n)
P_{n-1}(x;q;a,b,c).
\tag{5.3}
$$
More information, including the orthogonality relations, concerning the
Al-Salam--Chihara polynomials can be found in Askey and Ismail
\cite{5, \S 3}.

So we obtain the explicit expressions
$$
S_n(x) = \sum_{k=0}^n \frac{x^{n-k} q^{k(k-1)/2}}{(q;q)_k}
P_k(-(q+q^\nu);q;-(1+q^{\nu-1})q^{n-k+2},q^{2(n-k)+\nu+3},q^{\nu+1})
\tag{5.4}
$$
and
$$
\align
R_n(x) &=  \sum_{k=0}^n \frac{x^{n-k} q^{k(k-1)/2}}{(q;q)_k}
P_k(-(1+q^\nu);q;-(1+q^\nu)q^{n-k+1},q^{2(n-k)+\nu+2},q^\nu) \\
&\quad +\ \sum_{k=0}^{n-1} \frac{x^{n-1-k} q^{k(k+1)/2}}{(q;q)_k}
P_k(-(1+q^\nu);q;-(1+q^\nu)q^{n-k},q^{2(n-1-k)+\nu+2},q^\nu)   \\
&= x^n + \sum_{k=1}^n \frac{x^{n-k}q^{k(k-1)/2}}{(q;q)_k}
\left[ P_k(-(1+q^\nu);q;-(1+q^\nu)q^{n-k+1},q^{2(n-k)+\nu+2},q^\nu)
           \right. \\
&\qquad\qquad + (1-q^k) \left.
P_{k-1}(-(1+q^\nu);q;-(1+q^\nu)q^{n-k+1},q^{2(n-k)+\nu+2},q^\nu) \right].
\endalign
$$
A generating function for the Al-Salam--Chihara polynomials is,
cf. \cite{1, p.~23},
$$
\Phi(x,z) = \sum_{n=0}^\infty P_n(x;q;a,b,c) \frac{z^n}{(q;q)_n}
    = \frac{(\alpha z;q)_\infty (\beta z;q)_\infty}{(\gamma z;q)_\infty
        (\delta z;q)_\infty},
$$
where $1-az+bz^2=(1-\alpha z)(1-\beta z)$ and $1-xz+cz^2 =
(1-\gamma z)(1-\delta z)$.
Take $x=-(1+q^\nu)$ and $c=q^\nu$ so that
$\gamma=-1$ and $\delta = -q^\nu$. Consequently
$(1+z)\Phi(z,-(1+q^\nu))$ is the generating function for
$x=-(q+q^\nu)$ and $c=q^{\nu+1}$. Hence,
$$
P_n(-(1+q^\nu);q;a,b,q^\nu) + (1-q^n) P_{n-1}(-(1+q^\nu);q;a,b,q^\nu)
   = P_n(-(q+q^\nu);q;a,b,q^{\nu+1}),
$$
and thus
$$
R_n(x) = \sum_{k=0}^n \frac{x^{n-k} q^{k(k-1)/2}}{(q;q)_k}
     P_k(-(q+q^\nu);q;-(1+q^\nu)q^{n-k+1},q^{2(n-k)+\nu+2},q^{\nu+1}).
\tag{5.5}
$$

Now that we have the explicit expression for the polynomials $P_n$ defined in
\thetag{5.1} at hand, we can determine the asymptotic behaviour, which is
related to the function $j_\nu(x;q)$ introduced in \thetag{3.5}.

\proclaim{Proposition~5.1}
For the orthogonal polynomials $P_n(x)$ defined by \thetag{5.1}
we have for every $x\in\C$
$$
\lim_{n \to \infty} x^n P_n(1/x) = x^{1-\nu} j_{\nu-1}(x;q).
$$
\endproclaim

\demo{Proof} We follow the proof of theorem~2 of \cite{20}. For this we
need the continuous $q$-Hermite polynomials $H_n(x\mid q)$ introduced by
Rogers in 1894. The three-term recurrence relation is
$$
H_{n+1}(x\mid q) = 2x H_n(x\mid q) - (1-q^n) H_{n-1}(x\mid q)
\tag{5.6}
$$
with initial conditions $H_{-1}(x\mid q)=0$, $H_0(x\mid q)=1$, cf. Askey and
Ismail \cite{4, \S 6}. From \thetag{5.3} and \thetag{5.6} we
obtain, cf. \cite{20, thm.~2},
$$
\lim_{n\to\infty} P_k(-a;q;-aq^{n-k+1},b^2q^{2(n-k)+1},b^2/q) =
(-1)^k b^kq^{-{k\over 2}} H_k\Bigl( {{aq^{1/2}}\over{2b}}\mid q\Bigr).
\tag{5.7}
$$

Using this limit relation and dominated convergence we obtain
$$
\align
&\lim_{n \to \infty} x^n R_n(1/x) = \lim_{n \to \infty} x^n S_n(1/x) \\
&= \sum_{k=0}^\infty  \frac{(-1)^k q^{k(k-1)/2}}{(q;q)_k} x^k
q^{{k\over 2}(\nu+1)} H_k\Bigl({1\over 2}(
q^{{1\over 2}(1-\nu)}+q^{{1\over 2}(\nu-1)})\mid q\Bigr),
\endalign
$$
and hence
$$
\lim_{n \to \infty} x^n P_n(1/x) = \sum_{k=0}^\infty
\frac{(-1)^k q^{k(k-1)/2}}{(q;q)_k} x^{2k}
q^{{k\over 2}(\nu+1)} H_k\Bigl({1\over 2}(
q^{{1\over 2}(1-\nu)}+q^{{1\over 2}(\nu-1)})\mid q\Bigr).
\tag{5.8}
$$

To see that the right hand side of \thetag{5.8} equals
$x^{1-\nu}j_{\nu-1}(x;q)$ we insert the explicit expression,
cf. \cite{4, (6.1), (3.1)},
$$
H_k \Bigl( {1\over 2}(x+x^{-1})\mid q\Bigr) = \sum_{l=0}^k
{{(q;q)_k}\over{(q;q)_l(q;q)_{k-l} }} x^{k-2l}
$$
for $x=q^{(\nu-1)/2}$ in \thetag{5.8}. Interchanging summations and
introducing $m=k-l$ shows that \thetag{5.8} equals
$$
\align
&\sum_{l=0}^\infty {{(-1)^l q^{{1\over 2}l(l+1)}x^{2l}}\over{(q;q)_l}}
\sum_{m=0}^\infty {{(-1)^m q^{{1\over 2}m(m-1)}}\over{(q;q)_m}}
x^{2m} q^{m(l+\nu)} = \\
&\sum_{l=0}^\infty {{(-1)^l q^{{1\over 2}l(l+1)}x^{2l}}\over{(q;q)_l}}
(x^2q^{l+\nu};q)_\infty =
(x^2q^\nu;q)_\infty \, {}_1\vp_1 (0;x^2q^\nu;q,qx^2)
= x^{1-\nu}j_{\nu-1}(x;q)
\endalign
$$
by use of \cite{11, (1.3.16)}.\qed\enddemo

Observe that the continuous $q$-Hermite polynomials are orthogonal on the
interval $[-1,1]$, so that the inequality
$2\leq q^{{1\over 2}(1-\nu)}+q^{{1\over 2}(\nu-1)}$ shows that the
variable of the continuous $q$-Hermite polynomial in \thetag{5.8}
lies outside the support of the
orthogonality measure for the continuous $q$-Hermite polynomials, except
when $\nu=1$ in which case it is an endpoint of the interval.

The Stieltjes transform of the orthogonality measure $\mu$ for the orthogonal
polynomials $P_n$ can be obtained from
$$
\int_\R {{d\mu(t)}\over{z-t}} = \lim_{n\to\infty} {{P^{(1)}_{n-1}(z)}\over{P_n(z)}}
\tag{5.9}
$$
where $P^{(1)}_n$ are the associated polynomials, cf. Askey and Ismail
\cite{5, thm.~2.4} and further references therein.

So let us now consider the associated monic polynomials $P^{(1)}_n$ satisfying
$$
P^{(1)}_{n+1}(x) = x P^{(1)}_n(x) - \gamma_n P^{(1)}_{n-1}(x),
\qquad P^{(1)}_{-1}(x)=0, \quad P^{(1)}_0(x)=1,
\tag{5.10}
$$
where $\gamma_n=\lambda_{n+1}$ is defined in \thetag{5.1}.
These polynomials can be determined in a similar fashion as before.
Because of the parity of these polynomials, we again set
$$
P^{(1)}_{2n}(x) = T_n(x^2), \quad P^{(1)}_{2n+1}(x) = x U_n(x^2),
$$
and the monic polynomials $T_n$ and $U_n$ then satisfy the recurrence
relations
$$
\align
T_{n+1}(x) &= (x-\gamma_{2n}-\gamma_{2n+1})T_n(x) -
     \gamma_{2n-1}\gamma_{2n} T_{n-1}(x), \\
U_{n+1}(x) &= (x-\gamma_{2n+1}-\gamma_{2n+2})U_n(x) -
     \gamma_{2n}\gamma_{2n+1} U_{n-1}(x),
\endalign
$$
$T_0(x)=1$, $T_1(x)=x-q$ and $U_{-1}(x)=0$, $U_0(x)=1$, where
$$
\gamma_{2n}+\gamma_{2n+1} = \cases (q+q^\nu) q^n & \text{if $n>0$}, \\
                                              q & \text{if $n=0$.}
                                      \endcases , \quad
\gamma_{2n}\gamma_{2n-1} = q^{2n+\nu}, \qquad n \geq 0,
$$
and
$$
\gamma_{2n+1}+\gamma_{2n+2} = (1+q^\nu) q^{n+1}, \quad
\gamma_{2n}\gamma_{2n+1} = q^{2n+\nu+1}, \qquad n \geq 0.
$$

Hence
$$
\align
U_n(x) &= u_n(x;q(1+q^\nu),q^{(\nu+3)/2};q)  \\
      &=  \sum_{k=0}^n \frac{x^{n-k} q^{k(k-1)/2}}{(q;q)_k}
       P_k(-q(1+q^\nu);q;-(1+q^\nu)q^{n-k+2},q^{2(n-k)+\nu+4},q^{\nu+2}).
\endalign
$$
The polynomials $T_n$ are again co-recursive polynomials for the
recurrence relation \thetag{5.2} with $a=q+q^\nu$ and $b=q^{(\nu+2)/2}$, with
$T_1(x) = u_1(x) + q^\nu$, and thus
$$
T_n(x) = u_n(x;q+q^\nu,q^{(\nu+2)/2};q) +
 q^{\nu+n-1} u_{n-1}(x/q;q+q^\nu,q^{(\nu+2)/2};q).
$$
From the generating function of the Al-Salam--Chihara polynomials we find
$$
\multline
   P_n(-q(1+q^\nu);q;a,b,q^{\nu+2}) = \\
  P_n(-(q+q^\nu);q;a,b,q^{\nu+1}) + q^\nu (1-q^n)
  P_{n-1}(-(q+q^\nu);q;a,b,q^{\nu+1}) ,
  \endmultline
$$
so that
$$
T_n(x) =  \sum_{k=0}^n \frac{x^{n-k} q^{k(k-1)/2}}{(q;q)_k}
        P_k(-q(1+q^\nu);q;-(q+q^\nu)q^{n-k+1},q^{2(n-k)+\nu+3},q^{\nu+2}).
$$

The proof of the following proposition is analogous to the proof of
proposition~5.1.

\proclaim{Proposition~5.2}
For every $x\in\C$ we have
$$
\lim_{n \to \infty} x^n P^{(1)}_n(1/x) = x^{-\nu} j_\nu(x;q).
$$
\endproclaim

\subhead 5.2. Zeros of $j_\nu(x;q)$ and orthogonality relations\endsubhead
Combining propositions~5.1 and 5.2 and \thetag{5.9} shows that the Stieltjes
transform of the orthogonality measure $\mu$ for the polynomials is
$$
\int_\R {{d\mu(t)}\over{z-t}} = {{j_\nu(1/z;q)}\over{j_{\nu-1}(1/z;q)}}
\tag{5.11}
$$
for all $z\not\in {\hbox{\rm supp}}(d\mu)$, cf. \cite{5, thm.~2.4}.
From the Stieltjes transform we can derive the orthogonality relations
for the orthogonal polynomials $P_n$ defined in \thetag{5.1}. We start
with an investigation of the zeros of $j_\nu(x;q)$. It turns out that
the zeros of the function $j_\nu(x;q)$ behave similarly to the zeros of
the (Hahn-Exton $q$-)Bessel function for $\nu>-1$. The method of proof
largely follows Ismail's investigation \cite{15} of the roots of
the Jackson $q$-Bessel function, see also \cite{2, \S 4}.

\proclaim{Theorem~5.3} Let $\nu\in\R$ and let the function $j_\nu(x;q)$ be
defined by \thetag{3.5}.

\noindent
{\rm (i)} The functions $j_\nu(x;q)$ and $j_{\nu+1}(x;q)$ have no common
zeros, except possibly $x=0$.

\noindent
{\rm (ii)} The zeros of $x^{-\nu}j_\nu(x;q)$ are real, simple and symmetric
with respect to $x=0$. There are infinitely many of them and their only point
of accumulation is $\infty$.

\noindent
{\rm (iii)} The zeros of $x^{-\nu}j_\nu(x;q)$ and $x^{-\nu-1}j_{\nu+1}(x;q)$
interlace. Moreover, the smallest positive zero of $x^{-\nu}j_\nu(x;q)$ is
smaller than the smallest positive zero of $x^{-\nu-1}j_{\nu+1}(x;q)$.
\endproclaim

\demo{Proof} First we prove (i) by use of an equality for
${}_1\vp_1$-series. The relation
$$
{}_1\vp_1(0;c;q,z) - {}_1\vp_1(0;c;q,qz) = {{-z}\over{1-c}}\,
{}_1\vp_1 (0;cq;q,qz)
\tag{5.12}
$$
can be proved directly or it can be obtained from one of Heine's
contiguous relations for the
${}_2\vp_1$-series, cf. \cite{11, ex.~1.9(iv)}. Take $c=qx^2$ and
$z=q^{\nu+1} x^2$ in \thetag{5.12} to get from \thetag{3.5}
$$
j_\nu(x;q) -x^{-1} j_{\nu+1}(x;q) =
-q^{1+\nu/2} x^2 j_\nu(x\sqrt{q};q).
\tag{5.13}
$$
Substituting $c=q^{\nu+2}x^2$, $z=qx^2$ in \thetag{5.12} and using
\thetag{3.5} gives
$$
j_{\nu+1}(x;q) - q^{-{1\over 2}\nu} x j_\nu(x\sqrt{q};q) =
-q^{{1\over 2}(1-\nu)} x^2 j_{\nu+1}(x\sqrt{q};q).
\tag{5.14}
$$

If $0\not= a$ is a zero of $j_\nu(x;q)$ and $j_{\nu+1}(x;q)$, then
\thetag{5.13} implies that $a\sqrt{q}$ is a zero of $j_\nu(x;q)$. Next
\thetag{5.14} implies that $a\sqrt{q}$ is a zero of $j_{\nu+1}(x;q)$ as
well. So $aq^{k/2}$, $k\in\Zp$, are zeros of the analytic function
$x^{-\nu}j_\nu(x;q)$, which implies that this function is zero. This
contradiction proves (i).

To prove (ii) we recal that
the orthogonality measure $d\mu$ is supported on a
bounded denumerable discrete set
with zero as the only point of accumulation. So let $d\mu$ have
mass $A_k$ at the points $\{ t_k\}_{k=1}^\infty$, then \thetag{5.11} is
$$
\sum_{k=1}^\infty {{A_k}\over{z-t_k}} =
{{j_\nu(1/z;q)}\over{j_{\nu-1}(1/z;q)}}, \qquad z\not= t_k.
\tag{5.15}
$$
The zeros of $x^{1-\nu}j_{\nu-1}(1/x;q)$ correspond precisely to the
non-zero poles $t_k$ of the left hand side. So the zeros are real and simple.
Since $\{ t_k\}_{k=1}^\infty$ has zero as the only point of
accumulation, the only point of accumulation of the zeros of
$j_{\nu-1}(x;q)$ is infinity.

To prove (iii) we consider the (positive) mass of $d\mu$ at a non-zero $t_k$,
$$
0<A_k = -t_k^2 {{j_\nu(1/t_k;q)}\over{j_{\nu-1}^\prime (1/t_k;q)}}.
$$
So $j_\nu(a;q)$ and $j_{\nu-1}^\prime(a;q)$ have opposite signs for
$0\not=a$ a zero of $j_{\nu-1}(x;q)$. If $0<a<b$ are two consecutive
zeros of $j_{\nu-1}(x;q)$, then
$j_{\nu-1}^\prime(a;q)j_{\nu-1}^\prime(b;q)<0$. Hence also
$j_\nu(a;q)j_\nu(b;q)<0$ and $j_\nu(x;q)$ has at least one zero in
$(a,b)$. In the interval $(1/b,1/a)$ both sides of \thetag{5.15} are
differentiable, and the derivative of the left hand side is strictly
negative. If $j_\nu(1/z;q)$ has more than one zero in $(1/b,1/a)$, then
the derivative  has a zero in that interval. Thus $j_\nu(x;q)$ has
precisely one zero in $(a,b)$. This proves the interlacing property.

Denote by $x_k^\nu$ the positive zeros of $j_\nu(x;q)$ numbered
increasingly;
$$
0<x_1^\nu < x_2^\nu < \ldots < x_j^\nu < x_{j+1}^\nu < \ldots.
$$
Then it remains to prove that $x_1^{\nu-1}<x_1^\nu$. Since
$x^\nu j_\nu(x;q)$ equals $1$ for $x=0$ we get that
$j_{\nu-1}^\prime(x_1^{\nu-1};q)<0$ and thus
$j_\nu(x_1^{\nu-1};q)>0$. So $j_\nu(x;q)$ has an even number of zeros in
$(0,x_1^{\nu-1})$, and the same argument as in the previous paragraph
shows that this number is zero.
\qed\enddemo

The following proposition is the analogue of proposition~4.1 for the
orthogonal polynomials $P_n$ and the functions $j_\nu(x;q)$.

\proclaim{Proposition~5.4} For $n\in\Zp$ the polynomials $P_n$ and
$P^{(1)}_n$ defined by \thetag{5.1} and \thetag{5.10} satisfy
$$
P_n\bigl({1\over x}\bigr) j_\nu(x;q) -
P^{(1)}_{n-1}\bigl({1\over x}\bigr) j_{\nu -1}(x;q) =
\cases
{\displaystyle{ q^{m(m+{1\over 2}\nu)} x^{2m} j_\nu(xq^{{1\over 2}m};q), }}
                  &\text{$n=2m$,} \\
{\displaystyle{ q^{m(m+{1\over 2}(\nu-1))} x^{2m}
j_{\nu-1}(xq^{{1\over 2}m};q), }}
                  &\text{$n=2m-1$,}
\endcases
$$
where $j_\nu(x;q)$ is defined in \thetag{3.5}.
\endproclaim

\demo{Proof} It suffices to show that the right hand side satisfies
\thetag{5.1} with $x$ replaced by $x^{-1}$, since the left hand side
satisfies this equation and the cases $n=0$ (trivial) and $n=1$ (from
\thetag{5.13}) are easily proved. For $n=2m$ this follows from \thetag{5.13}
with $x$, $\nu$ replaced by $xq^{m/2}$, $\nu-1$, and for $n=2m-1$ this
follows from \thetag{5.14} with $x$, $\nu$ replaced by $xq^{(m-1)/2}$,
$\nu-1$.  \qed \enddemo

In the proof of theorem~5.3 we obtained information on the orthogonality
measure for the polynomials $P_n$ defined in \thetag{5.1}. In the next
theorem we describe the full orthogonality relations. This theorem can also
be proved from proposition~5.4 analogously to the proof of theorem~4.2 from
proposition~4.1.

\proclaim{Theorem~5.5} Let $\nu\in\R$, denote by $x_k^{\nu-1}$, $k\in\N$,
the positive zeros of the function $j_{\nu-1}(x;q)$ defined in
\thetag{3.5}. Then for the polynomials $P_n$ defined by \thetag{5.1} we
have the orthogonality relations
$$
\align
&\sum_{k=1}^\infty
P_n\Bigl( {{\pm 1}\over{x_k^{\nu-1}}}\Bigr)
P_m\Bigl( {{\pm 1}\over{x_k^{\nu-1}}}\Bigr)
{{-j_\nu(x_k^{\nu-1};q)}\over{(x_k^{\nu-1})^2
j_{\nu-1}^\prime(x_k^{\nu-1};q)}} \\
&\quad + (1-q^{\nu-1}) P_n(0)P_m(0) = \delta_{n,m}
\cases {\displaystyle{ q^{{1\over 2}l(l+\nu)},}}
        &\text{$n=2l$,} \\
       {\displaystyle{ q^{{1\over 2}(l+1)(l+\nu)},}}
        &\text{$n=2l+1$,}
\endcases
\endalign
$$
where the mass at $x=0$ only occurs for $\nu>1$. All weights are
positive.
\endproclaim

\demo{Proof} The only statements to be proved concern the norm and the
weight at $x=0$. Denote the squared norm of $P_n$ by
$\parallel P_n\parallel^2$, then \thetag{5.1} implies, cf.
\cite{10, (7)},
$$
\parallel P_n\parallel^2 = \lambda_n \parallel P_{n-1}\parallel^2 \qquad
\Longrightarrow \qquad
\parallel P_n\parallel^2 = \lambda_n\ldots\lambda_1
\parallel 1\parallel^2 .
$$
Together with the explicit value for $\lambda_n$ in \thetag{5.1} the
statement on the norm follows if we prove $\parallel 1\parallel^2=1$.
The value of $\parallel 1\parallel^2$ can be read off from the Stieltjes
transform \thetag{5.11} as the coefficient of $z^{-1}$ on the right hand
side and lemma~3.3 for $n=0$ and $x=z^{-1}$ shows that it equals $1$.

The weight at $x=0$ equals $\rho$, where $\rho^{-1} = \sum_{n=0}^\infty
\tilde P_n(0)^2$ and $\tilde P_n$ denote the orthonormal polynomials,
cf. \cite{5, thm.~2.8}.
From \thetag{5.1} we compute
$P_{2n+1}(0) = 0$, $P_{2n}(0) = (-1)^n q^{n\nu + n(n-1)/2}$
so that for the orthonormal polynomials we have
$$
\tilde P_{2n}(0) = {{P_{2n}(0)}\over{
\sqrt{\lambda_1\lambda_2\ldots\lambda_{2n}} }}
= (-1)^n q^{n(\nu-1)/2},
$$
and thus
$$
\sum_{n=0}^\infty \tilde P_n^2(0) = \sum_{n=0}^\infty q^{n(\nu-1)}
   = \cases   \infty & \text{if $\nu \leq 1$}, \\
              (1-q^{\nu-1})^{-1} & \text{if $\nu > 1$}, \endcases
$$
so that there is a mass $1-q^{\nu-1}$ at the origin whenever $\nu > 1$.
\qed\enddemo

Again, as in the proof of corollary~4.3, using \cite{21, (1.3) with
$n\to\infty$} shows that the orthogonality measure for the $P_n$ is
contained in $[-N,N]$ with $N\leq 1+q^{\nu/2}$. Shifting $\nu$ to $\nu+1$
we get the following corollary.

\proclaim{Corollary~5.6} The first positive zero $x_1^\nu$ of $j_\nu(x;q)$
satisfies $x_1^\nu\geq \bigl( 1+q^{(\nu+1)/2}\bigr)^{-1}$.
\endproclaim

\subhead 5.3. The case $\nu = 1/2$ \endsubhead
In the simple case $\nu=1/2$ we have
$\lambda_n = q^{n/2}$. For simplicity we take $p=q^{1/2}$ so that the
recurrence relation \thetag{5.1} can be rewritten as
$$
P_{n+1}(x) = x P_n(x) - p^n P_{n-1}(x).
\tag{5.16}
$$
We consider the generating function $G(z,x) = \sum_{n=0}^\infty P_n(x) z^n$.
Multiply \thetag{5.16} by $z^{n+1}$ and add all the terms from $n=0$ to
infinity, then we get
$$
G(z,x) - 1 = xz G(z,x) - z^2 p G(pz,x) \Longrightarrow
G(z,x) = \frac{1}{1-xz} - \frac{z^2 p}{1-xz} G(zp,x).
$$
Solving the $p$-difference equation with respect to the condition $G(0,x)=1$
gives by iteration
$$
G(z,x) = \sum_{k=0}^\infty \frac{(-1)^k z^{2k} p^{k^2}}{(zx;p)_{k+1}}.
\tag{5.17}
$$

We use the $p$-binomial theorem,
$$
\frac{1}{(zx;p)_{k+1}} = \sum_{n=0}^\infty
    \frac{(p^{k+1};p)_n}{(p;p)_n} (zx)^n,
$$
in \thetag{5.17}. Changing the summation index $n$ to $j-2k$ gives
$$
\align
G(z,x) & = \sum_{k=0}^\infty \sum_{j=2k}^\infty (-1)^k z^j x^{j-2k} p^{k^2}
    \frac{(p^{k+1};p)_{j-2k}}{(p;p)_{j-2k}} \\
& = \sum_{j=0}^\infty z^j \sum_{k=0}^{\lfloor j/2\rfloor}
     (-1)^k x^{j-2k} p^{k^2} \frac{(p^{k+1};p)_{j-2k}}{(p;p)_{j-2k}}.
\endalign
$$
Next identify the coefficient of $z^n$ and use
$(p^{k+1};p)_{j-2k} = (p;p)_{j-k}/(p;p)_k$ to find
$$
P_n(x) = \sum_{k=0}^{\lfloor n/2 \rfloor}
 (-1)^k x^{n-2k} p^{k^2} \frac{(p;p)_{n-k}}{(p;p)_k (p;p)_{n-2k}} .
\tag{5.18}
$$

These polynomials are a special case of orthogonal polynomials associated with
the Rogers-Ramanujan continued fraction; they correspond to the case
$a=0$, $b=p$ and $q=p$ in Al-Salam and Ismail
\cite{2} and \thetag{5.18} corresponds to
\cite{2, (3.7)}.
These polynomials are also the special case $u_n(x)$ in Van Assche
\cite{20} with $a=0$, $b=q$ and $q^2=p$, and \thetag{5.18} corresponds
to \cite{20, (2.7)} after observing that for the Al-Salam--Chihara
polynomials in \thetag{5.3} we have
$$
P_{2n+1}(0;q;0,b,c) = 0, \quad P_{2n}(0;q;0,b,c) =
(-1)^n c^n (\frac{b}{c};q^2)_n (q;q^2)_n .
$$

The associated polynomials $P_n^{(1)}$ satisfy the recurrence relation
$$
P_{n+1}^{(1)}(x) = xP_n^{(1)}(x) - p^{n+1} P_{n-1}^{(1)}(x),
\tag{5.19}
$$
with $P_{-1}^{(1)}=0$ and $P_0^{(1)}(x)=1$. Replace $x$ by $x/\sqrt{p}$
in \thetag{5.19}, then the monic polynomials $p^{n/2}P_n(x/\sqrt{p})$
satisfy the recurrence relation \thetag{5.19} so that $P_n^{(1)}(x)=
p^{n/2}P_n(x/\sqrt{p})$.

In case $\nu=1/2$ we have two different expressions for the same
polynomials. From \thetag{5.4} and \thetag{5.18} we obtain the following
summation formula for the Al-Salam--Chihara polynomials, $0\leq k\leq n$,
$$
P_k(-(q+q^{1\over2});q;-(1+q^{-{1\over 2}})q^{n-k+2},q^{2(n-k)+{7\over
2}}, q^{{3\over 2}}) =
{{(-q^{1\over 2})^k(q;q)_k (q^{1\over 2};q^{1\over 2})_{2n+1-k}}\over
{(q^{1\over 2};q^{1\over 2})_k (q^{1\over 2};q^{1\over
2})_{2n+1-2k}}}.
\tag{5.20}
$$
The Al-Salam--Chihara polynomials are expressible in a ${}_3\vp_2$-series as
proved by Askey and Ismail \cite{5, \S 3.8}. Explicitly, the following
connection between the original notation of Al-Salam and Chihara \cite{1}
and the notation of Askey and Ismail \cite{5} holds;
$$
\gathered
{{\alpha^{-k}}\over{(q;q)_k}} P_k(2\alpha
x;q; (\gamma+\delta)\alpha , \gamma\delta\alpha^2 , \alpha^2) =
S_k(x;\gamma,\delta\mid q)\\
= {{(\gamma\delta;q)_k}\over{(q;q)_k}}\gamma^{-k}
\, {}_3\vp_2 \left( {{q^{-k}, \gamma y, \gamma/y}\atop{\gamma\delta,\ 0}}
;q,q\right),
\endgathered
\tag{5.21}
$$
where $x=(y+y^{-1})/2$.

\proclaim{Corollary~5.7} The summation formula
$$
(c^2;q)_k \ {}_3\vp_2 \left( {{q^{-k},cq^{-1/2},c}\atop{c^2,\ 0}}
;q,q\right) = (cq^{-1/2})^k\, (-q^{1/2};q^{1/2})_k
(c;q^{1/2})_k
$$
holds for $k\in\Zp$ and $c\in\C$.
\endproclaim

\demo{Proof} In \thetag{5.20} we use \thetag{5.21} with the
parameters $\alpha=-q^{3\over 4}$, $x=(q^{1\over 4}+q^{-{1\over 4}})/2$,
$\gamma=q^{n-k+{3\over 4}}$, $\delta=q^{n-k+{5\over 4}}$ to get the result of
the proposition for $c=q^{n-k+1}$. Replace $n-k$ by $m$ in this result,
so that we have proved the corollary for $c=q^{m+1}$, $m\in\Zp$. Since both
sides are polynomial in $c$, the result follows for arbitrary values of $c$.
\qed\enddemo

\demo{Remark} Comparison of \thetag{5.5} with \thetag{5.18} instead of
\thetag{5.4} with \thetag{5.18} leads to the same corollary.
The same result is also obtained if we work out the different
expressions for the associated polynomials in case $\nu=1/2$. \enddemo

\proclaim{Proposition~5.8}
Consider the monic orthogonal polynomials given by \thetag{5.16} and the
associated polynomials given by \thetag{5.19}.
Then for every $x\in\C$ we have
$$
\lim_{n \to \infty} x^n P_n(1/x) = F(x), \qquad
\lim_{n\to\infty} x^n P^{(1)}_n(1/x) = F(x\sqrt{p})
$$
where
$$
F(x) = \sum_{k=0}^\infty \frac{(-1)^k x^{2k} p^{k^2}}{(p;p)_k}
   = {}_0\vp_1\left({{-}\atop{0}};p,-x^2 p \right).
$$
\endproclaim

\demo{Proof}
Straightforward by letting $n \to \infty$ in \thetag{5.18} after changing
$x$ to $1/x$ and multiplication by $x^n$. \qed
\enddemo

From proposition~5.8 and propositions~5.1 and 5.2 for $\nu=1/2$
we obtain the equalities
$$
{}_0\vp_1 (-;0;q^{1\over 2},-x^2q^{1\over 2}) = x^{1\over 2}j_{-1/2}(x;q),
\qquad
{}_0\vp_1(-;0;q^{1\over 2},-x^2q) = x^{-{1\over 2}}j_{1/2}(x;q),
\tag{5.22}
$$
which gives two transformations of a ${}_0\vp_1$-series of base $q^{1/2}=p$
in terms of ${}_1\vp_1$-series of base $q$.
In two special cases the left hand sides of \thetag{5.22} can be summed by the
Rogers-Ramanujan identities, cf. e.g. \cite{11, (2.7.3), (2.7.4)} and
this gives explicit values for $j_{-1/2}(x;q)$ for $x=\pm i, \pm iq^{1/4}$
and for $j_{1/2}(x;q)$ for $x=\pm i, \pm iq^{-1/4}$.

On comparing powers of $x$ in the first equality of
\thetag{5.22} and using the expansion \thetag{5.8} with $\nu=1/2$ for
$x^{1-\nu}j_{\nu-1}(x;q)$ we find for the continuous $q$-Hermite polynomial
the evaluation
$$
H_k \Bigl( {1\over 2}(q^{1\over 4}+q^{-{1\over 4}})\mid q\Bigr) = q^{-k/4}
(-q^{1/2};q^{1/2})_k.
$$
The same result is obtained if we compare powers of $x$ in the second equality
of \thetag{5.22}. Moreover, this result is a limit case of corollary~5.7, cf.
\thetag{5.7}.


\head 6. Orthogonality for the Laurent $q$-Lommel polynomials
\endhead

In this section we give a different form for the strong moment functional
introduced in \S 3.
The limit transitions \thetag{3.4} and \thetag{3.6} suggest to rewrite the
strong moment functional $\L$ defined in \thetag{3.15} as a contour integral
over the unit circle. This can be done if $j_{\nu-1}(1;q)\not= 0$,
since we have sufficient knowledge
on the location of the zeros of $j_{\nu-1}(x;q)$, cf. theorem~5.3, and of
$J_{\nu-1}(x;q)$, cf. \cite{17, \S 3} and \S 4. A Wronskian type formula
can be used to simplify the integrand.

\proclaim{Lemma~6.1} Let $r_m(x)$, $s_m(x)$ be solutions of the recurrence
relation \thetag{1.7}, then the Wronskian
$r_m(x)s_{m+1}(x)-s_m(x)r_{m+1}(x)$ is independent of $m\in\Z$.
\endproclaim

\demo{Proof} Multiply the recurrence formula for $r_m(x)$ by $s_m(x)$ and
multiply the recurrence relation for $s_m(x)$ by $r_m(x)$. Subtract the
resulting identities to find the result.
\qed\enddemo

\proclaim{Lemma~6.2}
$$
J_\nu(1/x;q)j_{\nu-1}(x;q) - J_{\nu-1}(1/x;q)j_\nu(x;q) = x^{-1}
{{(qx^{-2};q)_\infty (x^2;q)_\infty}\over{(q;q)_\infty}}.
$$
\endproclaim

\demo{Proof}
$J_{\nu+m}(1/x;q)j_{\nu+m-1}(x;q) - J_{\nu+m-1}(1/x;q)j_{\nu+m}(x;q)$
is independent of $m$ by pro\-po\-si\-tion 3.1 and lemma~6.1.
Take $m=0$ to obtain the left hand side of the lemma
and use \thetag{3.7} and $m\to\infty$ to see that it also equals
$$
(x^{-1}-x){{(qx^{-2};q)_\infty (qx^2;q)_\infty}\over{(q;q)_\infty}},
$$
which proves the lemma.
\qed\enddemo

Lemma~6.2 implies that
$x(q;q)_\infty \bigl[ J_\nu(1/x;q)j_{\nu-1}(x;q) -
J_{\nu-1}(1/x;q)j_\nu(x;q)\bigr]$
is a theta product, cf. Askey \cite{3, \S 1}.

Now we can rewrite the strong moment functional $\L$
with respect to which the Laurent
$q$-Lommel polynomials are orthogonal, cf. theorem~3.4.

\proclaim{Theorem~6.3} Let $s>0$ such that $s$ is not a zero of
$J_{\nu-1}(1/x;q)$ and $j_{\nu-1}(x;q)$.  For $\nu>0$ the strong moment
functional $\L$ defined in \thetag{3.15} equals
$$
\gather
\L (p) = {1\over{2\pi i(q;q)_\infty}} \oint_{\vert z\vert = s}
p(z)\, {{(qz^{-2};q)_\infty (z^2;q)_\infty}\over{
J_{\nu-1}(1/z;q) j_{\nu-1}(z;q)}} \, {{dz}\over z}\ \ + \\
\sum_{k=1}^N \Bigl( p\bigl( {1\over{j_k^{\nu-1}}}\bigr)
+ p\bigl( {{-1}\over{j_k^{\nu-1}}}\bigr)\Bigr)
{{-J_\nu(j_k^{\nu-1};q)}\over{(j_k^{\nu-1})^2
J_{\nu-1}^\prime(j_k^{\nu-1};q)}}
+ \sum_{l=1}^M \bigl( p( x_l^{\nu-1})
+ p( -x_l^{\nu-1})\bigr)
{{j_\nu(x_l^{\nu-1};q)}\over{ j_{\nu-1}^\prime(x_l^{\nu-1};q)}},
\endgather
$$
where $p$ is an arbitrary Laurent polynomial. Here $j_k^{\nu-1}$, respectively
$x_l^{\nu-1}$, denote the positive zeros of $J_{\nu-1}(x;q)$, respectively
$j_{\nu-1}(x;q)$, numbered increasingly. $N$ is defined by
$j_N^{\nu-1}<s<j_{N+1}^{\nu-1}$ and $N=0$, and so the sum over the zeros
of $J_{\nu-1}(x;q)$ is empty, if $j_1^{\nu-1}> s$.
$M$ is defined by
$x_M^{\nu-1}<s< x_{M+1}^{\nu-1}$ and $M=0$, and so the sum over the zeros
of $j_{\nu-1}(x;q)$ is empty, if $x_1^{\nu-1}> s$. The discrete weights
in the first sum over $k$ are positive and the discrete weights in the second
sum over $l$ are negative.
\endproclaim

\demo{Remark} (i) By choosing $s=r$, respectively $s=1/R$, with $r$, $R$ as
in \S 3, we get $M=0$, respectively $N=0$. In \S 7 we show that for $\nu$
sufficiently large we have $N=M=0$ for a suitable choice of $s$.

\noindent
(ii) The non-zero poles of the integrand in theorem~6.3 are simple. Indeed,
if $0\not= a$ satisfies $J_{\nu-1}(1/a;q)=0=j_{\nu-1}(a;q)$, then lemma~6.2
implies that the numerator is zero as well. Moreover, $a=q^{p/2}$ for some
$p\in\Z$, which is a simple zero of the numerator. There exist only finitely
many of such values in the (possibly empty) interval
$[x_1^{\nu-1},1/j_1^{\nu-1}]$.
\enddemo

\demo{Proof} In the first contour integral in \thetag{3.15} we shift the
contour integration from $\vert z\vert = 1/R$ to $\vert z\vert =s$
and in general we assume
$s<1/R$. We pick up residues at the simple poles
$z=\pm 1/j_k^{\nu-1}$, $k=1,\ldots,N$,
cf. theorem~4.2. For $1/R\leq s$ we are in the case $N=0$.
The second contour integral in \thetag{3.15} is shifted from $\vert z\vert =r$
to $\vert z\vert =s$. In general we assume $r<s$, otherwise we are in case
$M=0$. Here we pick up residues at the simple poles
$z=\pm x_l^{\nu-1}$, $l=1,\ldots,M$. The residues are easily calculated.
Next we take together the integrands of the contour integrals over
$\vert z\vert = s$ using lemma~6.2 to prove the expression for $\L(p)$ in
this case.
The last statement follows from theorem~4.2 and theorem~5.5.
\qed\enddemo

\demo{Remark} The most natural choice for $s$ in theorem~6.3 seems $s=1$.
This is motivated by the fact that on the unit circle there is a
transition in
the asymptotic behaviour of the Laurent $q$-Lommel polynomials,
cf. \thetag{3.4}, \thetag{3.6}. Moreover, numerical experiments
indicate that for $m\to\infty$
the non-real zeros of the Laurent $q$-Lommel polynomials, cf. remark~2.1,
are possibly dense on the unit circle. Of course, from \thetag{3.4},
respectively \thetag{3.6}, we see that the real zeros outside, respectively
inside, the unit circle tend to the zeros of $J_{\nu-1}(x^{-1};q)$,
respectively $j_{\nu-1}(x;q)$. This corresponds precisely with the discrete
set in the orthogonality measure of theorem~6.3 for $s=1$.
\enddemo


\head 7. Laurent $q$-Lommel polynomials as perturbations of Chebyshev
polynomials
\endhead

Let us now return to the recurrence relation \thetag{1.7} which we rewrite as
$$
h_{n+1,\nu}(x;q) - (x^{-1} + x) h_{n,\nu}(x;q) + h_{n-1,\nu}(x;q) =
     -xq^{\nu+n} h_{n,\nu}(x;q) .
\tag{7.1}
$$
In this way, as $q \to 0$ or as $\nu \to \infty$ the Laurent polynomials
$h_{n,\nu}(x;q)$  should be close to a solution of the three term recurrence
relation
$$
h_{n+1}(x;0) - (x^{-1} + x) h_{n}(x;0) + h_{n-1}(x;0) = 0.
\tag{7.2}
$$
The solution of this recurrence, with initial values $h_0(x;0)=1$ and
$h_{-1}(x;0)=0$ is given by $h_n(x;0) = (x^{n+1}-x^{-n-1})/(x-x^{-1})$, which
in terms of Chevbyshev polynomials of the second kind can be written as
$$
h_n(x;0) = U_n\x , \qquad n \in \Zp.
$$
In this way the Laurent polynomials $h_{n,\nu}(x;q)$ can be considered
as perturbations of the Chebyshev polynomials. We now do a perturbation
analysis, much as is done for perturbations of orthogonal polynomials in
\cite{21}. In the spirit of the
Liouville-Green approximation (WKB method), we will consider \thetag{7.1}
as a second order recurrence relation with non-homogeneous term
$-xq^{\nu+n} h_{n,\nu}(x;q)$, even though this term  depends on the
desired solution $h_{n,\nu}(x;q)$.

We solve this non-homogeneous recurrence relation by
Green's method. We need the Green function $G_1(n,m)$, which is the solution
of the recurrence relation with non-homogeneous term $\delta_{n,m}$, i.e.,
$$
G_1(n+1,m) - (x^{-1}+x) G_1(n,m) + G_1(n-1,m) = \delta_{n,m}
\tag{7.3}
$$
with boundary conditions
$$
G_1(n,m) = 0, \qquad n \geq m .
\tag{7.4}
$$
Clearly $G_1(m,m) = G_1(m+1,m) = 0$ and thus from \thetag{7.3} we
find $G_1(m-1,m) = 1$. For $k \geq 0$ we find that $r_k(x)=G_1(m-k-1,m)$ is a
solution of the homogeneous recurrence relation \thetag{7.2} with the
same initial conditions $r_0(x)=1$ and $r_{-1}(x)=0$, hence
$$
G_1(n,m) = U_{m-n-1} \x , \qquad n < m.
$$
Now multiply \thetag{7.1} by $G_1(n,m)$ and \thetag{7.3} by $h_{n,\nu}(x;q)$
and subtract the obtained equations to find
$$
\multline
   h_{n+1,\nu}(x;q) G_1(n,m) - h_{n,\nu}(x;q)G_1(n-1,m)
  +  h_{n,\nu}(x;q) \delta_{n,m}  \\
 =  h_{n,\nu}(x;q)G_1(n+1,m) - h_{n-1,\nu}(x;q)G_1(n,m)
     - xq^{\nu+n} h_{n,\nu}(x;q) G_1(n,m).
\endmultline
$$
Add all the equations from $n=0$ to $n=m$ and use the boundary conditions
\thetag{7.4} to find
$$
h_{0,\nu}(x;q)G_1(-1,m) =   h_{m,\nu}(x;q) + x \sum_{n=0}^{m-1}
    q^{\nu+n} G_1(n,m) h_{n,\nu}(x;q).
$$
This gives
$$
h_{m,\nu}(x;q) = U_m \x - x \sum_{n=0}^{m-1} q^{\nu+n} U_{m-n-1} \x
   h_{n,\nu}(x;q) .
\tag{7.5}
$$
From this relation we can deduce some useful properties.

\proclaim{Lemma~7.1}
Suppose $x = e^{i\theta}$ with $\theta \in [0,2\pi)$, then
$$
|h_{n,\nu}(x;q)| \leq (n+1) \exp \left( \frac{q^\nu}{(1-q)^2} \right),
\tag{7.6}
$$
and
$$
|\sin \theta\ h_{n,\nu}(x;q) | \leq 1 +  \frac{q^\nu}{(1-q)^2}
  \exp \left( \frac{q^\nu}{(1-q)^2} \right) .
\tag{7.7}
$$
For $|x| \neq 1$ we have
$$
\aligned
| x^n h_{n,\nu}(x;q) | &\leq \frac{2}{|1-x^2|} \exp \left(
    \frac{2}{|1-x^2|} \frac{q^\nu}{1-q} \right), \qquad |x| < 1, \\
  | x^{-n} h_{n,\nu}(x;q) &| \leq \frac{2}{|1-x^{-2}|} \exp \left(
    \frac{2}{|1-x^{-2}|} \frac{q^\nu}{1-q} \right)
  , \qquad |x| > 1.
\endaligned
$$
\endproclaim

\demo{Proof}
We use Gronwall's inequality, cf. \thetag{3.9}; for
non-negative $A, c_n, d_n$, $(n \geq 0)$ we have
$$
c_n \leq A + \sum_{k=0}^{n-1} d_k c_k \Longrightarrow
c_n \leq A \exp \left( \sum_{k=0}^{n-1} d_k \right).
$$
From the bound $|U_n(\cos \theta)| \leq n+1$ and \thetag{7.5} we find
$$
|h_{n,\nu}(x;q)| \leq n+1 + \sum_{k=0}^{n-1} q^{\nu+k} (n-k)
     |h_{k,\nu}(x;q)|.
$$
Hence taking $c_n = |h_{n,\nu}(x;q)|/(n+1)$ in Gronwall's inequality
gives
$$
\frac{|h_{n,\nu}(x;q)|}{n+1} \leq \exp \left( \sum_{k=0}^{n-1} (k+1)
    q^{\nu +k} \right).
$$
The desired inequality \thetag{7.6}
then follows from $\sum_{k=0}^\infty (k+1)q^k =
(1-q)^{-2}$. If we use this inequality \thetag{7.6} and
$|\sin \theta\ U_n(\cos \theta)| \leq 1$ in \thetag{7.5}, then
$$
|\sin \theta\ h_{n,\nu}(x;q)| \leq 1 + \exp \left( \frac{q^\nu}{(1-q)^2}
    \right) \sum_{k=0}^{n-1} (k+1)   q^{\nu+k} ,
$$
which gives \thetag{7.7}.
The bounds away from the unit circle follow by using
$$
| x^n U_n \x| = \left| x^n \frac{x^{n+1}-x^{-n-1}}{x-x^{-1}} \right|
\leq \frac{2}{|1-x^2|}, \qquad |x| < 1,
$$
and
$$
| x^{-n} U_n \x| = \left| x^{-n} \frac{x^{n+1}-x^{-n-1}}{x-x^{-1}} \right|
\leq \frac{2}{|1-x^{-2}|}, \qquad |x| > 1,
$$
and by using Gronwall's inequality. \qed
\enddemo

From these bounds we see that the Laurent Lommel polynomials have an
exponentially increasing upper bound both inside the unit circle and
outside the unit circle, and that on the unit circle the Laurent
polynomials are bounded,
except when $x=\pm 1$, in which case $|h_{n,\nu}(x;q)| = O(n)$.
This strongly suggests that in theorem~6.3 the choice $s=1$ for
the strong moment functional $\L$ is the most natural.

The Laurent polynomial solution of \thetag{7.1} is not the only interesting
solution. In \S 3 we already obtained the minimal solutions
$j_{\nu+n}(x;q)$ and $J_{\nu+n}(x^{-1};q)$ on respectively
the open unit disk and the exterior of the closed unit disk. The minimal
solutions $h_n^-(x;0)$ and $h_n^+(x;0)$ of the recurrence relation
\thetag{7.2} on respectively the open unit disk $\{ z \in \C:\ |z| < 1 \}$
and the exterior of the closed unit disk $\{ z \in \C:\ |z| > 1 \}$ are given
by $h_n^-(x;0) = x^n$ and $h_n^+(x;0) = x^{-n}$. Our intention now is
to find similar solutions $h_{n,\nu}^\pm(x;q)$ satisfying
$$
\lim_{n \to \infty} h_{n,\nu}^\pm(x;q) x^{\pm n} = 1 ,
$$
on $\{ z \in \C:\ |z| < 1 \}$ and $\{ z \in \C:\ |z| > 1 \}$ respectively.
Such functions clearly exist, since by proposition~3.1 and \thetag{3.7}
we see that
$$
h_{n,\nu}^+(x;q) = \frac{(q;q)_\infty}{(qx^{-2};q)_\infty}
x^{\nu}J_{\nu+n}(x^{-1};q) , \quad
h_{n,\nu}^-(x;q) = \frac{1}{(qx^2;q)_\infty}
x^{-\nu}j_{\nu+n}(x;q) ,
\tag{7.8}
$$
fulfil the required conditions.

We will now do a perturbation analysis of these minimal solutions
in a similar way as is done for orthogonal polynomials \cite{12}.
Again we write the recurrence relation as
$$
h_{n+1,\nu}^\pm - (x^{-1} + x) h_{n,\nu}^\pm(x;q) + h_{n-1,\nu}^\pm(x;q) =
 -xq^{\nu+n} h_{n,\nu}^\pm(x;q) ,
\tag{7.9}
$$
and look at this equation as a non-homogeneous second order
recurrence relation
with non-homogeneous term $-xq^{\nu+n} h_{n,\nu}^\pm(x;q)$.
The homogeneous equation
has two simple solutions, $h_n^\pm(x;0) =  x^{\mp n}$.
We solve the non-homogeneous recurrence relation using Green functions, but
now the Green function $G_2(n,m)$ is the solution of
$$
G_2(n+1,m) - (x^{-1}+x) G_2(n,m) + G_2(n-1,m) = \delta_{n,m}
\tag{7.10}
$$
with boundary conditions
$$
G_2(n,m) = 0, \qquad n \leq m.
\tag{7.11}
$$
Since $G_2(m,m)=G_2(m-1,m)=0$ we find $G_2(m+1,m)=1$ and in general
$$
G_2(n,m) = U_{n-m-1} \x, \qquad n > m.
$$

Multiply the recurrence \thetag{7.9} by $G_2(n,m)$ and
\thetag{7.10} by $h_{n,\nu}^\pm(x;q)$ and subtract to find
$$
\multline
h_{n+1,\nu}^\pm(x;q) G_2(n,m) - h_{n,\nu}^\pm(x;q)G_2(n-1,m)
+ h_{n,\nu}^\pm(x;q) \delta_{n,m} \\
= h_{n,\nu}^\pm(x;q)G_2(n+1,m) - h_{n-1,\nu}^\pm(x;q)G_2(n,m)
 - xq^{\nu+n} h_{n,\nu}^\pm(x;q) G_2(n,m).
\endmultline
$$
Add the equations from $n=m$ to $n=M$, with $m < M$ and use the boundary
conditions \thetag{7.11} to find
$$
\multline
h_{M+1,\nu}^\pm(x;q) G_2(M,m) - h_{M,\nu}^\pm(x;q)G_2(M+1,m) \\
= - h_{m,\nu}^\pm(x;q) - x \sum_{n=m+1}^M q^{\nu+n}
h_{n,\nu}^\pm(x;q) G_2(n,m).
\endmultline
$$
From \thetag{7.8} together with \thetag{3.7} we obtain
$$
\lim_{M \to \infty} h_{M+1,\nu}^+(x;q) G_2(M,m) - h_{M,\nu}^+(x;q)G_2(M+1,m)
= -x^{-m}, \qquad |x| > 1,
$$
and
$$
\lim_{M \to \infty} h_{M+1,\nu}^-(x;q) G_2(M,m) - h_{M,\nu}^-(x;q)G_2(M+1,m)
 -x^{m}, \qquad |x| < 1,
$$
so by letting $M \to \infty$ we have
$$
h_{n,\nu}^\pm(x;q) = x^{\mp n} -  x \sum_{k=n+1}^{\infty} q^{\nu+k}
h_{k,\nu}^\pm(x;q) U_{k-n-1} \x .
\tag{7.12}
$$
Compare these relations to \thetag{7.5}. We can find appropriate bounds
on these solutions and from this we can obtain bounds for
$J_{\nu+n}(x^{-1};q)$ and $j_{\nu+n}(x;q)$.

\proclaim{Lemma~7.2} If $x \neq \pm 1$ then
$$
\aligned
|x^n h_{n,\nu}^+(x;q)| &\leq \exp \left( \frac{2}{|1-x^{-2}|}
\frac{q^{\nu+n+1}}{1-q} \right), \qquad |x| \geq 1, \\
|x^{-n} h_{n,\nu}^-(x;q)| &\leq \exp \left( \frac{2}{|1-x^{2}|}
\frac{q^{\nu+n+1}}{1-q} \right), \qquad |x| \leq 1,
\endaligned
\tag{7.13}
$$
and
$$
\aligned
|x^nh_{n,\nu}^+(x;q)| \leq \exp \left( \frac{n q^{\nu+n+1}}{1-q}
+ \frac{q^{\nu+n+1}}{(1-q)^2} \right), \qquad |x| \geq 1, \\
|x^{-n}h_{n,\nu}^-(x)| \leq \exp \left( \frac{n q^{\nu+n+1}}{1-q}
+ \frac{q^{\nu+n+1}}{(1-q)^2} \right), \qquad |x| \leq 1.
\endaligned
\tag{7.14}
$$
\endproclaim

\demo{Proof} We now use a backward version of Gronwall's inequality:
for non-negative $A, c_n, d_n$ $(n \geq 0)$ we have
$$
c_n \leq A + \sum_{k=n+1}^\infty d_k c_k < \infty \Longrightarrow
   c_n \leq A \exp \left( \sum_{k=n+1}^\infty d_k \right).
$$
The inequalities \thetag{7.13} then follow  from
\thetag{7.12} and the inequalities
$|x^{\pm n} U_n \xs | \leq 2/|1-x^{\pm 2}|$,
which hold for $|x| \leq 1$ (for the $+$ sign) and $|x| \geq 1$
(for the $-$ sign).

Inequality \thetag{7.14} uses the inequality $|x^{\pm n}U_n\xs|
\leq n+1$ on $|x| \leq 1$ and $|x| \geq 1$ respectively. So from
\thetag{7.12} we get
$$
\multline
\vert x^n h^+_{n,\nu}(x;q)\vert \leq 1 + \sum_{k=n+1}^\infty
q^{\nu+k} \vert x^k h^+_{k,\nu}(x;q)\vert (k-n) \\
\leq 1 + \sum_{k=n+1}^\infty
k\, q^{\nu+k} \vert x^k h^+_{k,\nu}(x;q)\vert, \qquad \vert x\vert \geq 1,
\endmultline
$$
from which the first inequality of \thetag{7.14} follows by Gronwall's
inequality.
\qed \enddemo

We are now ready to give some information about the zeros of the functions
$h_{n,\nu}^{\pm}(x;q)$ inside and outside the open unit disk.

\proclaim{Theorem~7.3} The zeros of $h_{n,\nu}^{\pm}(x;q)$ are
all real. The function $h_{n,\nu}^+(x;q)$ has no zeros
in $\{ x \in \C:\  |x| \geq 1\}$ and $h_{n,\nu}^-(x;q)$ has no
zeros inside $\{ x \in \C:\  |x| \leq 1\}$, whenever $n \geq M(\nu,q)$, where
$$
M(\nu,q) =  -\nu-1 + 2 \frac{\ln(1-q)}{\ln q} -
\frac{1}{\ln q} .
\tag{7.15}
$$
In particular $h_{-1,\nu}^+(x;q)$ has at most $2M(\nu,q)+2$ zeros
in $\{ x \in \C:\  |x| \geq 1\}$ and $h_{-1,\nu}^-(x;q)$ has at most
$2M(\nu,q)+2$ zeros in $\{ x \in \C:\  |x| \leq 1\}$.
\endproclaim

\demo{Proof}
The reality of the zeros follows from the explicit representation
\thetag{7.8} and the reality of the zeros of the Hahn-Exton Bessel function
\cite{17, \S 3} and the zeros of $j_\nu(x;q)$, cf. theorem~5.3.
For an upper bound on the number of zeros, we use \thetag{7.12} to find
$$
1- x^{\pm n} h_{n,\nu}^\pm(x;q)
= \sum_{k=n+1}^\infty q^{\nu+k} x^{\pm k} h^\pm_{k,\nu}(x;q)
x^{\pm n \mp k +1} U_{k-n-1}\x .
$$
Use the inequality \thetag{7.14} and $|x^{\pm n} U_n\xs| \leq n+1$ to find
for $|x| \geq 1$
$$
| 1- x^{n} h_{n,\nu}^+(x;q) | \leq
\sum_{k=n+1}^\infty (k-n) q^{\nu+k}
\exp \left( \frac{k q^{\nu+k+1}}{1-q}
+ \frac{q^{\nu+k+1}}{(1-q)^2} \right)  ,
$$
and similarly for $|x| \leq 1$
$$
| 1- x^{-n} h_{n,\nu}^-(x;q) | \leq
\sum_{k=n+1}^\infty (k-n) q^{\nu+k}
\exp \left( \frac{k q^{\nu+k+1}}{1-q}
 + \frac{q^{\nu+k+1}}{(1-q)^2} \right)  .
$$
The right hand side can be bounded by
$$
\align
\sum_{k=n+1}^\infty (k-n) q^{\nu+k}
   \exp \left( \frac{k q^{\nu+k+1}}{1-q}
 + \frac{q^{\nu+k+1}}{(1-q)^2} \right)
&\leq \exp \left(  \frac{q^{\nu+1}}{(1-q)^2} \right)
\sum_{k=n+1}^\infty (k-n) q^{\nu+k}  \\
&= \frac{q^{\nu+n+1}}{(1-q)^2} \exp \left(
\frac{q^{\nu+n+1}}{(1-q)^2} \right)  .
\endalign
$$
Choose $M=M(\nu,q)$ such that
$$
\frac{q^{\nu+M+1}}{(1-q)^2} \exp \left(
  \frac{q^{\nu+M+1}}{(1-q)^2} \right) < 1,
$$
then $h_{n,\nu}^+(x;q)$ for $n \geq M$ cannot be zero for any
$x$ such that $|x| \geq 1$. An appropriate $M(\nu,q)$ is given by
\thetag{7.15}. The same reasoning holds for $h_{n,\nu}^-(x;q)$ on
the closed unit disk.
So now we have established that for $n \geq M$ the function $h_{n,\nu}^+$
has no zeros for $|x| \geq 1$ and $h_{n,\nu}^-$ has no zeros for
$|x| \leq 1$. The zeros of $h_{n,\nu}^+$ are equal to the zeros
of $J_{\nu+n}(1/x;q)$. If $j_{k}^{\nu}$, $k=1,2,3,\ldots$, are
the zeros of $J_{\nu}(x;q)$ numbered increasingly, then from
the interlacing property of
theorem~3.7 in \cite{17} we have
$j_k^{\nu} < j_k^{\nu+1} < j_{k+1}^{\nu}$,
hence when the parameter $\nu$ is decreased by one, then the $k$th positive
zero moves to the left. This means that the $k$th positive
zero (counted from the right)
of $h_{n-1,\nu}^+(x;q)$ is to the right of the $k$th positive zero
of $h_{n,\nu}^+(x;q)$. Since $h_{M,\nu}^+(x;q)$ has no zeros  $x\geq 1$,
this means that $h_{M-1,\nu}^+(x;q)$ can have one zero $x \geq 1$, namely
$1/j_{1}^{\nu+M-1}$, and it
cannot have two zeros $x >1$ since $1/j_2^{\nu+M-1} < 1/j_1^{\nu+M} < 1$.
Decreasing the degree of $h_{n,\nu}^+(x;q)$ by one thus increases the number
of zeros in $|x| \geq 1$ by at most 2 (one positive zero and
one negative zero). Therefore $h_{-1,\nu}^+(x;q)$ has at most
$2M+2$ zeros in $|x| \geq 1$.
A similar reasoning works for the zeros of $h_{n,\nu}^-(x;q)$ in
$|x| \leq 1$ by using the interlacing property of the zeros of
$j_{\nu}(x;q)$ and $j_{\nu+1}(x;q)$ given by theorem~5.3.
\qed \enddemo

The upper bound on the number of zeros of $h_{-1,\nu}^\pm(x;q)$
gives a useful upper bound on the number of discrete mass points
of the strong moment functional $\L$ as given in theorem~6.3 when $s=1$.
Indeed, the zeros of $h_{-1,\nu}^+(x;q)$ correspond with the zeros of
$J_{\nu-1}(1/x;q)$ and thus $N \leq M(\nu,q)+1$. Similarly the zeros
of $h_{-1,\nu}^-(x;q)$ correspond with the zeros of $j_{\nu-1}(x;q)$ and thus
$M \leq M(\nu,q)+1$. In particular, $M=N=0$ in theorem~6.3 for $s=1$ for
$\nu$ satisfying $M(\nu,q)<0$.

Finally let us give another derivation of the orthogonality of the
Laurent polynomials $h_{n,\nu}(x;q)$ by using the minimal
solutions $h_{n,\nu}^\pm(x;q)$.
Observe that from \thetag{3.3}, \thetag{3.5} and \thetag{7.8}
it follows that
$h_{n,\nu}^\pm(x;q)$ have a power expansion of the form
$$
x^n h_{n,\nu}^+(x;q) = 1 + \sum_{k=0}^\infty K^+(n,k) x^{-2k},
 \qquad |x| > 1,
$$
and
$$
x^{-n} h_{n,\nu}^-(x;q) = 1 + \sum_{k=1}^\infty K^-(n,k) x^{2k},
  \qquad |x| < 1.
$$
We can get some information on the coefficients $K^\pm(n,k)$ by introducing
Banach algebras. If $f$ is analytic in the open unit disk with Taylor
series $$   f(z) = \sum_{k=0}^\infty f_k z^k, $$
then we define
$$
\|f \|_- = \sum_{k=0}^\infty \nu_k |f_k|,
$$
and we denote by $A^-$ all the functions $f$ for which $\| f\|_- < \infty$.
Here $\nu_k$, $k\in\Zp$, is a positive increasing sequence for which
$\nu_0=1$ and $\nu_{n} \leq \nu_m \nu_{n-m}$ for every $n \geq m \geq 0$.
Similarly, when $g$ is analytic near infinity with Laurent series
$$
g(z) = \sum_{k=0}^\infty  g_k z^{-k},
$$
then we define
$$
\|g \|_+ = \sum_{k=0}^\infty \nu_k |g_k|,
$$
and denote by $A^+$ all the functions $g$ for which $\| g\|_+ < \infty$.
One easily verifies that for two functions $f_1, f_2 \in A^\pm$ one has
$$
\| f_1 f_2 \|_\pm \leq \| f_1 \|_\pm\ \| f_2 \|_\pm,
$$
so that we are dealing with Banach algebras.

Observe that
$$
\gather
\| x^n h_{n,\nu}^+(x;q) \|_+ =
  1 + \sum_{k=0}^\infty \nu_{2k} |K^+(n,k)|,     \\
\| x^{-n} h_{n,\nu}^-(x;q) \|_- =
  1 + \sum_{k=1}^\infty \nu_{2k} |K^-(n,k)|.
\endgather
$$
Taking norms in \thetag{7.12} gives
$$
\| x^n h_{n,\nu}^+(x;q) \|_+ \leq
1 + \sum_{k=n+1}^\infty q^{\nu+k} \|x^k h_{k,\nu}^+(x;q) \|_+ \
\| x^{n-k+1} U_{k-n-1}\x \|_+ .
$$
Now
$$
\| x^{-n} U_n\x \|_+  = \| \sum_{j=0}^n x^{-2j} \|_+
= \sum_{j=0}^n \nu_{2j} \leq (n+1) \nu_{2n},
$$
so that Gronwall's inequality gives
$$
\| x^n h_{n,\nu}^+(x;q) \|_+ \leq
\exp \left( \sum_{k=n+1}^\infty
k q^{\nu+k}  \nu_{2k} \right).
$$
Taking $\nu_n = a^n$ with $a < q^{-1/2}$ shows that
$$
1 + \sum_{k=0}^\infty a^{2k} |K^+(n,k)| < \infty,
$$
so that $x^n h_{n,\nu}^+(x;q) \in A^+$. This shows that
the function $h_{n,\nu}^+(x;q)$ is in fact defined for
$x > q^{1/2}$. A similar reasoning shows that
$x^{-n} h_{n,\nu}^-(x;q) \in A^-$ and that $h_{n,\nu}^-(x;q)$ is defined
for $x < q^{-1/2}$. From \thetag{7.8} we see that
$h_{n,\nu}^+(x;q)$ has poles at the zeros of $(qx^{-2};q)_\infty$
and that $x=\pm q^{1/2}$ are the poles of largest modulus. Similarly
$h_{n,\nu}^-(x;q)$ has poles at the zeros of $(qx^2;q)_\infty$ and
$x=\pm q^{-1/2}$ are the poles of smallest modulus.

Suppose now that $\pm 1$ are not zeros of $h^\pm_{-1}(x;q)$.
Evaluate the contour integral
$$
I_+ = \frac{1}{2\pi i} \int_{|x|=1} \frac{h_{m,\nu}(x;q)
  h_{n,\nu}^+(x;q)}{h_{-1,\nu}^+(x;q)} \, dx  .
$$
If $m < n$ then near $x=\infty$ the integrand behaves as $x^{m-n-1}$ and thus
$I_+$ has no contribution from $x=\infty$. So
when $x_j^+$ $(j \geq 1)$ are the zeros of $h_{-1,\nu}^+(x;q)$, then
$$
I_+ = - \sum_{j=1}^{N} \frac{h_{m,\nu}(x_j^+;q)
h_{n,\nu}^+(x_j^+;q)}{[h_{-1,\nu}^+(x_j^+;q)]'} ,
$$
where $N$ is defined as in theorem~6.3 for $s=1$.
Similarly we compute the contour integral
$$
I_- = \frac{1}{2\pi i} \int_{|x|=1} \frac{h_{m,\nu}(x;q)
  h_{n,\nu}^-(x;q)}{h_{-1,\nu}^-(x;q)} \, dx  .
$$
The integrand behaves as $x^{n-m+1}$ near $x=0$ and thus there is no pole
at the origin when $m < n$ (even for $m \leq n+1$). There are poles at
the zeros $x_j^-$ $(j\geq 1)$ of $h_{-1,\nu}^-(x;q)$ and we thus have
$$
I_- = \sum_{j=1}^{M} \frac{h_{m,\nu}(x_j^-;q)
h_{n,\nu}^-(x_j^-;q)}{[h_{-1,\nu}^-(x_j^-;q)]'},
$$
where $M$ is defined as in theorem~6.3 for $s=1$.
Subtracting $I_+$ and $I_-$ gives
$$
I_- - I_+ = \frac{1}{2\pi i} \int_{|x|=1}
h_{m,\nu}(x;q) \left( \frac{ h_{n,\nu}^-(x;q)h_{-1,\nu}^+(x;q)
 - h_{n,\nu}^+(x;q)h_{-1,\nu}^-(x;q)}{h_{-1,\nu}^+(x;q)h_{-1,\nu}^-(x;q)}
\right) \, dx .
$$
The Laurent polynomial $h_{n,\nu}(x;q)$ is a solution of the three-term
recurrence relation \thetag{7.1} and therefore a linear combination of the
two special solutions $h_{n,\nu}^\pm(x;q)$. With the initial conditions
$h_{0,\nu}(x;q) = 1$ and $h_{-1;\nu}(x;q)=0$ and by combining
\thetag{7.8} with lemma~6.2 we find
$$
(x^{-1}-x) h_{n,\nu}(x;q) = h_{-1,\nu}^-(x;q) h_{n,\nu}^+(x;q)
- h_{-1,\nu}^+(x;q) h_{n,\nu}^-(x;q),
$$
so that
$$
I_- - I_+ = \frac{1}{2\pi i} \int_{|x|=1}
(x-x^{-1}) h_{n,\nu}(x;q) h_{m,\nu}(x;q) \frac{dx}{h_{-1,\nu}^+(x;q)
h_{-1,\nu}^-(x;q)} .
$$
On the other hand, at a zero $x_j^+$ we see that $h_{n,\nu}^+(x_j^+;q)$
is a solution of \thetag{7.1} with initial value $h_{-1,\nu}^+(x_j^+;q)=0$,
so that $h_{n,\nu}^+(x_j^+;q) = h_{0,\nu}^+(x_j^+;q) h_{n,\nu}(x_j^+;q)$.
Similarly at a zero $x_j^-$ we have
$h_{n,\nu}^-(x_j^-;q)=h_{0,\nu}^-(x_j^-;q) h_{n,\nu}(x_j^-;q)$. Therefore
$$
\align
I_- - I_+ &= \sum_{j=0}^{M(\nu,q)}
h_{n,\nu}(x_j^-;q) h_{m,\nu}(x_j^-;q) \frac{h_{0,\nu}^+(x_j^-;q)}
{[h_{-1,\nu}^+(x_j^-;q)]'}    \\
&\quad + \sum_{j=0}^{M(\nu,q)}
h_{n,\nu}(x_j^+;q) h_{m,\nu}(x_j^+;q) \frac{h_{0,\nu}^+(x_j^+;q)}
{[h_{-1,\nu}^+(x_j^-;q)]'}  .
\endalign
$$
Combining both expressions for $I_--I_+$ gives the orthogonality of the
Laurent polynomials $h_{n,\nu}(x;q)$ and corresponds to the result
given in theorem~6.3 for $s=1$. In case $q^{1/2}<s<q^{-1/2}$ the
orthogonality relations of theorem~6.3 can be derived in a similar way.

This approach can also be used to prove the orthogonality for the
Laurent polynomials $x^{-1}h_{n,\nu}(x;q)$, cf. \thetag{2.3}. Note also
that the case $q=0$ gives the orthogonality relations for the
Chebyshev polynomials of the second kind.


\enddocument